\documentstyle[12pt,graphicx,amsfonts]{article}
\newcommand{\vs}{\vspace{0.3cm}}

\newtheorem{theorem}{Theorem}[section]

\def\thebibliography#1
{\begin{center}{\normalsize\bf References}\end{center}
\list{[\arabic{enumi}]}%
{\settowidth\labelwidth{[#1]}\leftmargin \labelwidth
\advance\leftmargin\labelsep
\usecounter{enumi}}
\sloppy\clubpenalty4000\widowpenalty4000}
\begin{document}
\baselineskip=18pt

\title{\vspace*{-2cm}\large SIGNATURE of ROTORS}

\author{{\normalsize MIECZYS{\L}AW~ D\c{A}BKOWSKI}\\
{\small Department of Mathematics, University of Texas at Dallas}\\[-1mm]
{\small Richardson, Texas 75083-0688, USA}\\[-1mm]
{\small e-mail: mdab@utdallas.edu}
\\[3mm]
{\normalsize MAKIKO ISHIWATA}\\
{\small Department of Mathematics, Tokyo Woman's Christian University}\\[-1mm]
{\small Zempukuji 2-6-1, Suginamiku, Tokyo 167-8585, Japan}\\[-1mm]
{\small e-mail: mako@lab.twcu.ac.jp}
\\[3mm]
{\normalsize J\'{O}ZEF H. PRZYTYCKI}\\
{\small Department of Mathematics, The George Washington University}\\[-1mm]
{\small Washington, DC 20052, USA}\\[-1mm]
{\small e-mail: przytyck@gwu.edu}
\\[3mm]
{\normalsize AKIRA YASUHARA }\\
{\small Department of Mathematics, Tokyo Gakugei University}\\[-1mm]
{\small Nukuikita 4-1-1, Koganei, Tokyo 184-8501, Japan}\\[-1mm]
{\small e-mail: yasuhara@u-gakugei.ac.jp}
}

%\date{}
 \date{\today}
%\date{July 10, 2004}
\maketitle

\begin{abstract}
Rotors were introduced as a generalization of mutation 
by Anstee, Przytycki and Rolfsen in 1987. 
In this paper we show that Tristram-Levine signature is preserved by 
 orientation-preserving rotations.
Moreover, we show that any link invariant obtained from the
 characteristic polynomial of Goeritz matrix, including Murasugi
 signature, is not changed by rotations.
In 2001, P. Traczyk showed that the Conway polynomials of any pair of
orientation-preserving rotants coincide.
But it was still an open problem if an orientation-reversing rotation 
preserves Conway polynomial.
We show that there is a pair of orientation-reversing rotants 
with different Conway polynomials. 
This provides a negative solution to the problem.
\end{abstract}



{\section{Introduction}} 
{\it{Rotors}} were introduced in graph theory by 
W.Tutte \cite{BSST}, \cite{Tu1} and \cite{Tu2}. 
The concept was adapted to knot theory in \cite{APR} as a generalization of
 Conway's mutation. For the orientation of the  boundary of an 
{\it oriented rotor}, we have two basic possibilities:\\ 
(a) inputs and outputs alternate as 
in Fig.2.2 (a). Such a rotor is called an {\it orientation-preserving 
rotor}, or\\
 (b) we have the pattern in-in, out-out as in Fig.2.2 (b). 
Such a rotor is called an {\it orientation-reversing rotor}\footnote{The 
terminology used in here is explained in Section 2.}.
\par  
In Section 3 (resp. Section 4), we show, in particular, that the Murasugi's
unoriented version of the classical signature \cite{G-L, M, M2} (Theorem 
3.1) (resp. Tristram-Levine signature) is preserved by any rotations (resp. any orientation-preserving rotations).
\par
It was shown in \cite{APR} that rotations of order 
three and four preserve the Homflypt 
polynomial, and in particular, the Conway polynomial of links.
In 2001, P.Traczyk \cite{T} showed that Conway polynomials of a pair of
 any orientation-preserving rotants coincide, solving in 
this case, the Jin-Rolfsen Conjecture \cite{J-R}. 
But it was inconclusive if orientation-reversing 
rotations preserve Conway polynomials for $n \geq 6$. 
In the last section, we present an example of orientation-reversing 
rotants which do not share the same Conway 
polynomial. This provides a negative answer for the 
Jin-Rolfsen Conjecture in the orientation-reversing case \cite{J-R,P}.
%{\footnote{See Definition 2.3, {\textup
%{\cite{G-L, M, M2}}
%}.}} 
%???Furthermore we consider the first homology groups of the
%$p$-fold branched covers branched along any pair of 
%oriented rotants of any order $n$.
%\par
%What we know about preserving homology of double cover? 
%\par
%
%Remark:
\par
In general, it is not true that a rotation preserves 
the first homology of
the double branched cover, ${M_L}^{(2)}$, of $S^3$ branched along
$L$.%$(S^3,L)$.
~Necessary conditions for preserving the homology are given 
in \cite{DJP,P2}. 
%Memo: determinant: same, so the order of $H_1$.
%\par
Figure 1.1 taken from \cite{DJP} shows rotants with different
 $H_1(M_{L_k}^{(2)};\Bbb Z$) and $H_1(M^{(2)}_{L_k};{\Bbb Z}_5)$. 
%and $H_1(M_{L_k}^{(2)};\Bbb Z$) $(k=1,2)$.  
%{\footnote{$H_1(M_{L_k}^{(2)};\Bbb Z)$,$(k=1,2)$ are computed by 
%KODAMA software.}} 
%and $H_1(M^{(2)}_{L},\Bbb Z$$_5)$.
For the link $L_1$ in
Fig 1.1(a), $H_1(M_{L_1}^{(2)};\Bbb Z$$)=\Bbb Z$$_{15} \oplus \Bbb Z$$_{30}$ 
and $H_1(M_{L_1}^{(2)};\Bbb Z$$_5)=\Bbb Z$$_5 \oplus \Bbb Z$$_5$, and 
for its orientation preserving rotant
$L_2$ in Fig.1.1(b) we obtain
 $H_1(M_{L_2}^{(2)};\Bbb Z$$)=\Bbb Z$$_3 \oplus \Bbb 
Z$$_{150}$, $H_1(M_{L_2}^{(2)};\Bbb Z$$_5)=\Bbb Z$$_5$.
All the homology groups were calculated using K. Kodama's program 
KNOT \cite{K}. 
%KODAMA software.
\par\vs
\begin{center}
\begin{tabular}{cc} 
\includegraphics[trim=0mm 0mm 0mm 0mm, width=.25\linewidth]
{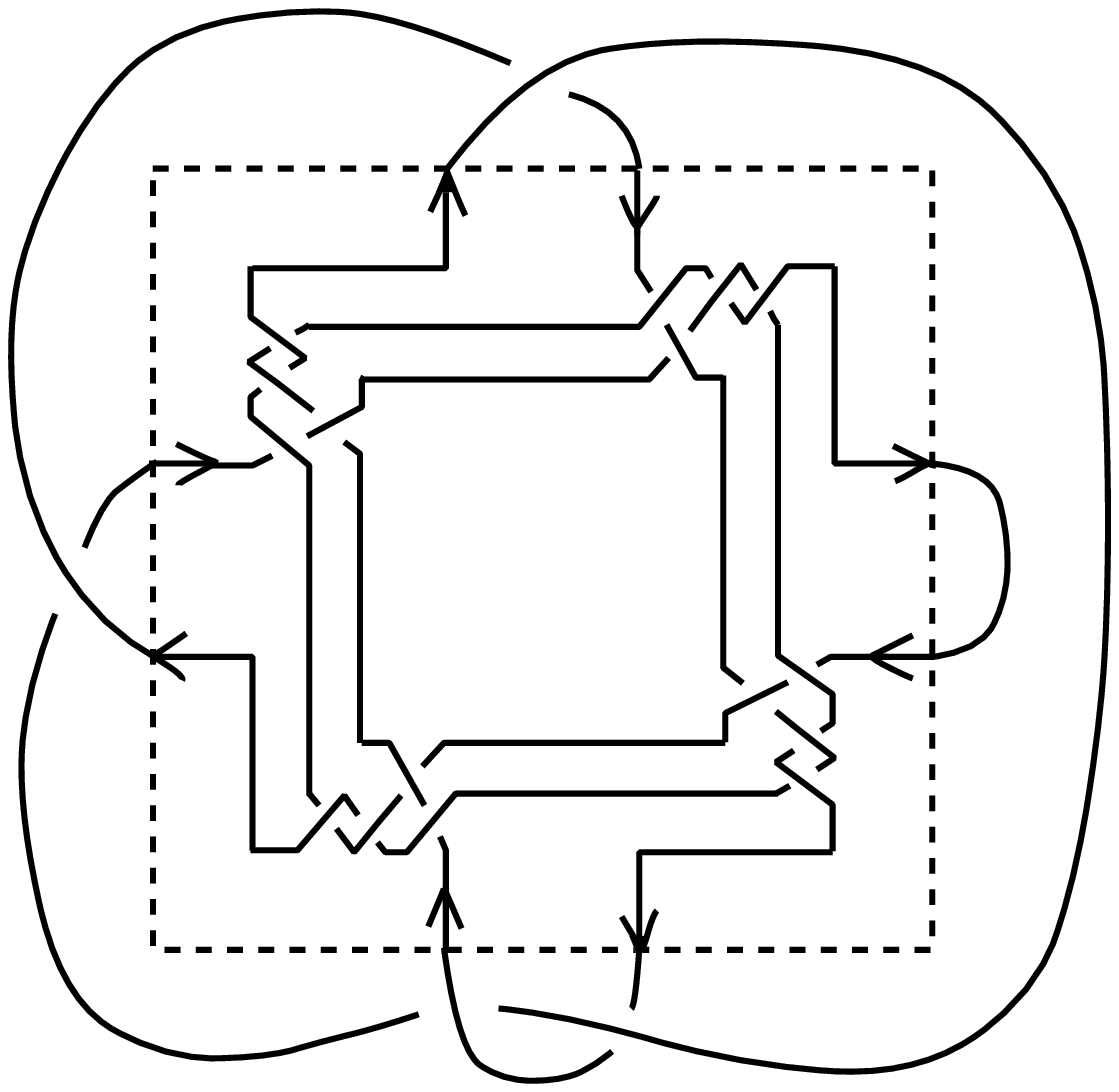}
&\hspace*{10mm}
\includegraphics[trim=0mm 0mm 0mm 0mm, width=.25\linewidth]
{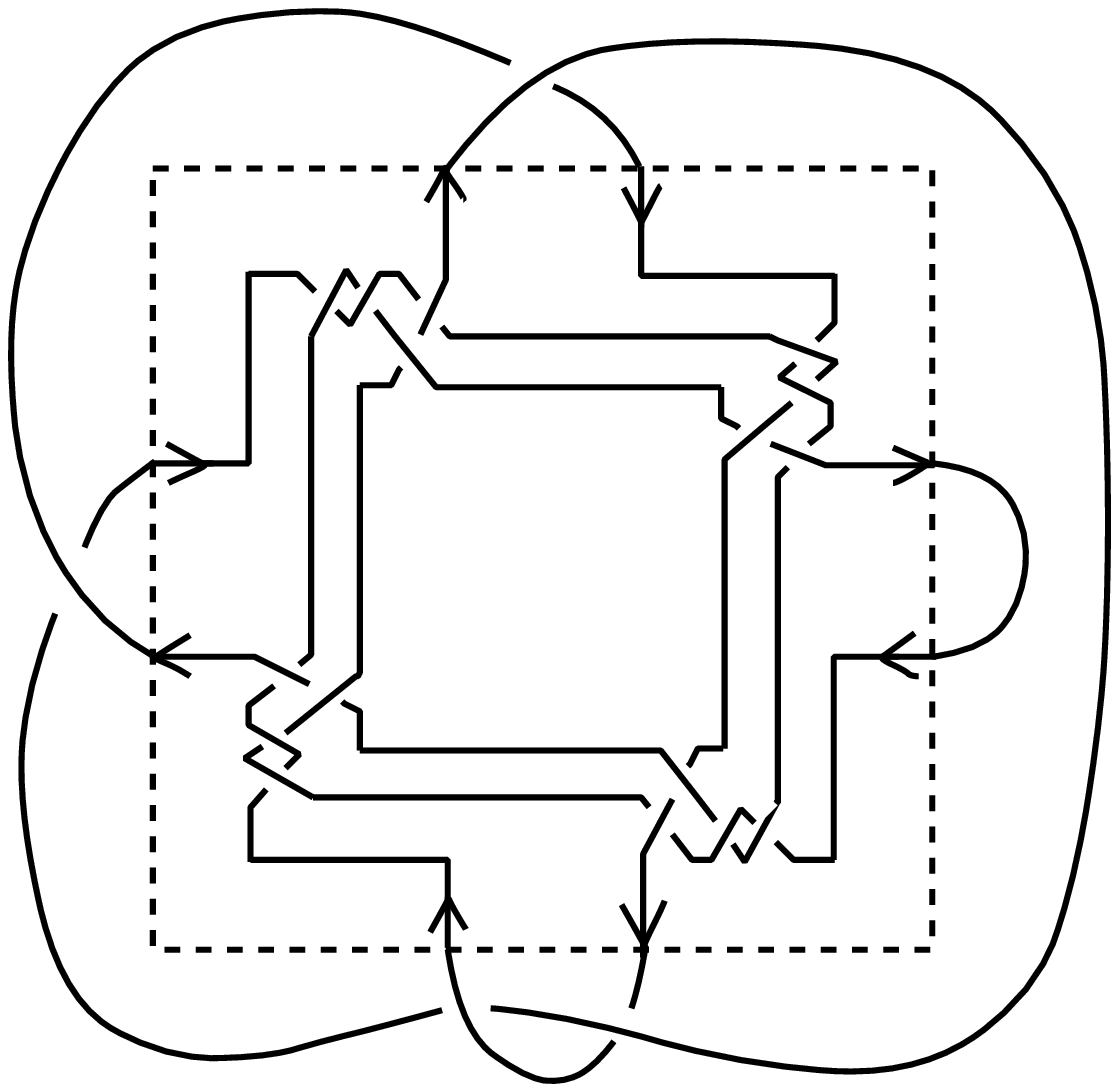}\\
$L_1$&\hspace*{10mm}$L_2$\\
(a)&\hspace*{10mm}(b)
\end{tabular}
\par
Fig. 1.1
\end{center}
However, if we assume that a given  pair of oriented rotants can 
be put into the ``special'' periodic 
disk-band form then the first homology groups of the corresponding
 double branched covers of $S^3$ branched along this pair of rotants 
%link $L_k (k=1,2)$ 
are isomorphic (Corollary 2.3).

\medskip
\noindent
{\section{Definitions and basic properties of rotors}}
\par
%%%%%%%%%%%%Definitions and basic properties%%%%%%%%%%%%%%%%%%%%%%%%
%The linking number is the oldest invariant of links. 
For an oriented
link $L$ of $k$-components $K_1,\cdots,K_k$ we form the linking matrix
$A_L$ with entries $a_{ij}=lk(K_i,K_j)$, where $i \ne j$. 
We put $a_{ii}=0$ unless $L$ is a framed link.
 In this case we define $a_{ii}$ to be the framing of the $i$th
component $K_i$ of $L$ ($a_{ii}$ measures the difference 
with respect to the standard framing).  The linking
%symmetric 
matrix $A_L$, up to the order of
components of $L$, is a link invariant. 
%The dimension of $A_L$ is the number of components of $L$.
 One half of the sum ${\sum_{i < j}}a_{ij}$ of entries of $A_L$
outside the diagonal is the total linking number of $L$, denoted by 
${\ell} k(L)$. The trace of $A_L$ for a framed link $L$ is denoted by
$tr(L)$. Note that $tr(L)$ does not depend on the orientation of $L$, so 
$tr(L)$ is an invariant of an unoriented framed link $L$.
\par
Consider a link $L$ in $S^3$ decomposed into two $n$-tangles ($n > 2$)
 $S$ and $R$ ~(Fig. 2.1), where by $n$-tangle we mean  any
 1-dimensional manifold 
properly embedded into a three-ball and consisting of $n$-arcs and, 
possibly, closed components.
Let $\phi$ be a rotation of $B^3=B^2 \times I$ by the 
angle $\frac{2 \pi}{n}$ along the $z$ axis.
Assume that $R$, called the {\it{rotor part}} of $L$, satisfies
$\phi(R)=R$.
% has a rotational symmetry. 
The other tangle part, $S$, of $L$ is called the {\it{stator}}. 
%Such a link 
Equivalently, $L$ admits a projection decomposed into the projections  
%$S$ and $R$  
of the rotor and the stator (these projections will also be
denoted by $S$ and $R$) such that $R$ lies in the regular $n$-gon and  
intersects 
its boundary in $2n$ points, and that 
% If $\phi$ is a rotation of the projection
%plane by the angle $\frac{2 \pi}{n}$, then 
$\phi(R)=R$ (Fig. 2.1).
\par
\vs\vs
The regular $n$-gon has a dihedral group of symmetry $D_{2n}$. This
group is generated by the $2\pi/n$ rotation along the $z$ axis  $\phi$ 
and the dihedral flype $d_0$  which corresponds to
the rotation by $\pi$ along the $y$ axis. 
The group $D_{2n}$ has a presentation, $D_{2n}=\{\phi, d_0
\mid {\phi}^n={d_0}^2=1, d_0 \phi d_0={\phi}^{-1}\}$. Let
$d_{\frac{k}{2}}={\phi}^k d_0$. 
Note that $d_{\frac{k}{2}}$ is the dihedral flype along the axis
obtained from the $y$ axis by rotating it counterclockwise
 by the angle $\frac{2\pi k}{2n}$. 
%the angle $\frac{2\pi k}{2n}$.
\par
A {\it{rotant}} of a link $L_1$ is the link $L_2$ (Fig. 1.1 and Fig. 2.1)
obtained from $L_1$ by a dihedral flype of its rotor part. Note that
$L_2$ is independent of the choice of a dihedral flype. We say that $L_2$ is
obtained from $L_1$ by a {\it{rotation}}.
\par\vs\vs\vs
\begin{center}
\begin{tabular}{cc} 
\includegraphics[trim=0mm 0mm 0mm 0mm, width=.30\linewidth]
{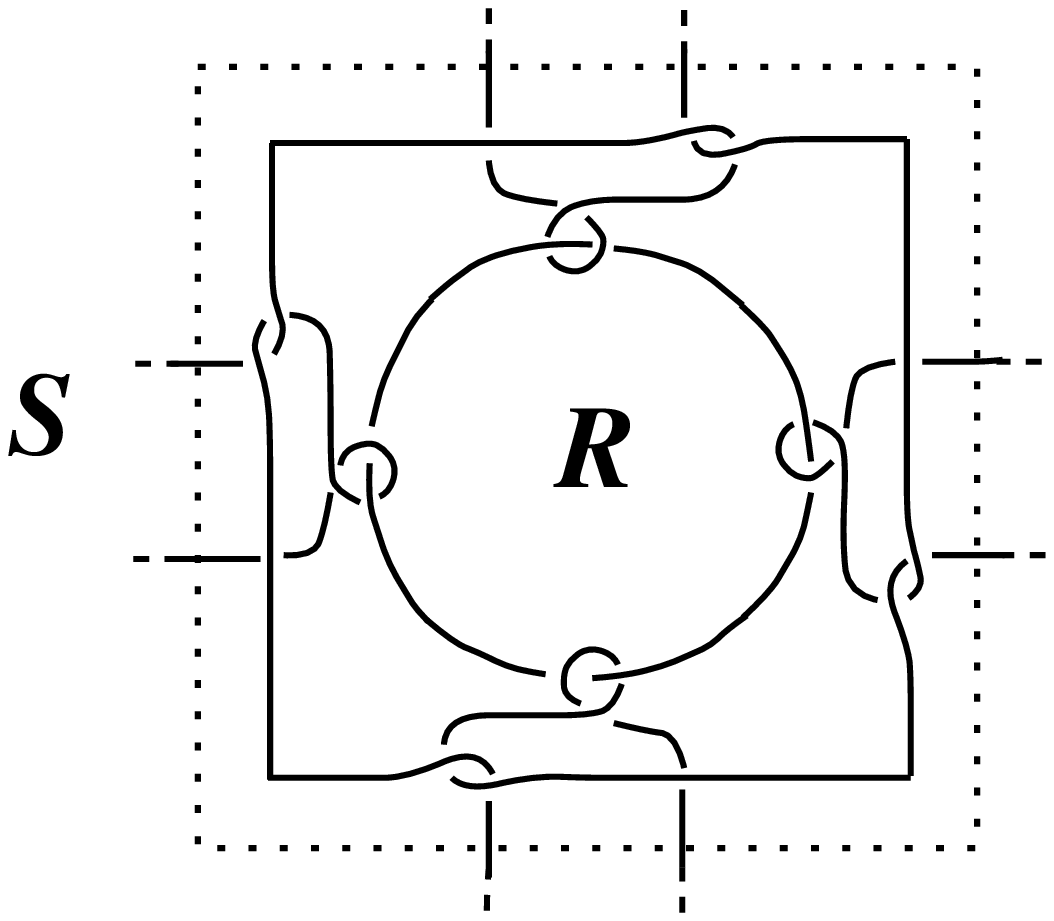}
&\hspace*{20mm}
\includegraphics[trim=0mm 0mm 0mm 0mm, width=.30\linewidth]
{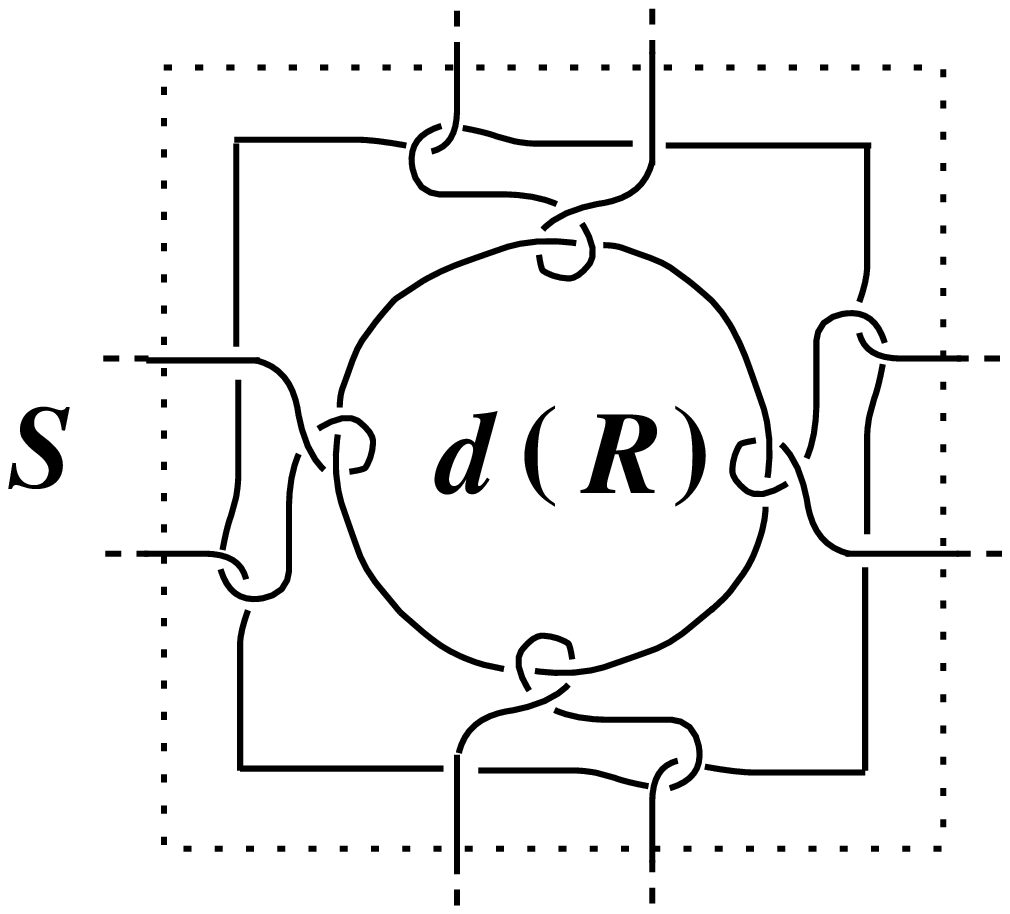}\\
$L_1$&\hspace*{20mm}$L_2$
\end{tabular}
\\
Fig. 2.1
\end{center}
\par\vs\vs
If a link is equipped with additional structures such as orientation or a 
blackboard framing, we also
assume that the rotation preserves these structures. 
In the oriented case, we allow the global change of the 
orientation of the rotor part. 
More precisely, for an oriented rotor we have two basic 
choices of directions of arcs at its boundary points: 
inputs and outputs alternate as in 
Fig. 2.2(a), we call such a rotor the {\it orientation-preserving rotor}, 
or we have the pattern in-in, out-out, $\cdots$, in-in, out-out 
for an even $n$ as in Fig. 2.2(b); we call such a rotor 
 the {\it{orientation-reversing rotor}}.  
For an oriented rotor $R$ of an oriented link $L$ and a dihedral flype $d$, 
the orientations of $d(R)$ and the stator parts do not  always 
necessarily match. If they do not match, then by reversing the 
orientation of $d(R)$, we obtain an
oriented link $L_2=d(R) \cup S$ that we also call the oriented 
rotant of $L_1$.
\par\vs
\begin{center}
\begin{tabular}{cc} 
\includegraphics[trim=0mm 0mm 0mm 0mm, width=.2\linewidth]
{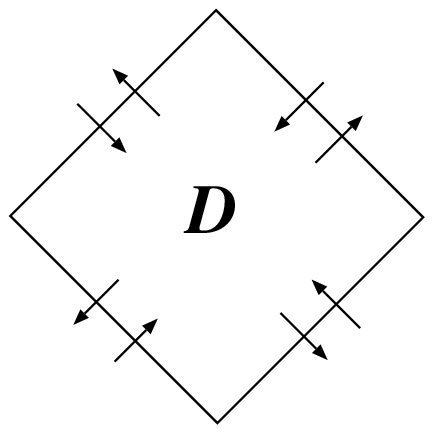}
&\hspace*{20mm}
\includegraphics[trim=0mm 0mm 0mm 0mm, width=.2\linewidth]
{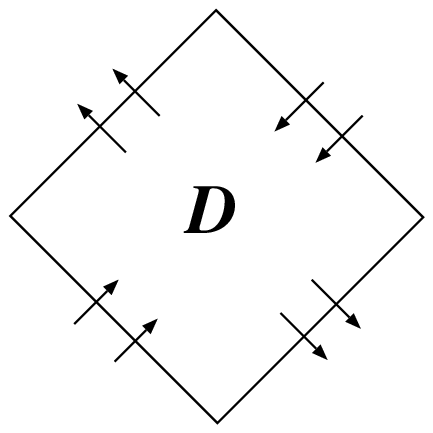}
\end{tabular}
\begin{tabular}{cccc}
(a)& \hspace*{45mm} (b)  
\end{tabular}\\
Fig. 2.2
\end{center}
\par\vs
The following theorem describes basic properties of rotors. 
\begin{theorem}
\begin{enumerate}
\item[{\em(i)}] Any rotation preserves the number of components of a link.

\item[{\em(ii)}] If two oriented links are related by a rotation of an	
		oriented rotor, then the total linking numbers are the same.

\item[{\em (iii)}] If two oriented framed links are related by a
		rotation of an oriented rotor, and the rotor part has no closed components,
		then their linking matrices are the same.

\item[{\em(iv)}] If $L$ is an unoriented framed link, then $tr(L)$ is
		preserved by any rotation.

%\item[\em{(v)}] If $L$ is an unoriented link diagram 
%then any rotation
%preserves the self-twist number $sw(L)$ of $L$, where $sw(L)$ 
%is defined to be the sum of signs of selfcrossings of $L$.

\end{enumerate}
\end{theorem}
\par\vs\vs
\par\vs\noindent
{\bf{Proof}}~~ 
%\par
Let $R$ be an unoriented rotor with boundary points $a_0,b_0,a_1,b_1,...,
a_{n-1}, b_{n-1}$, as in Figure 2.3(a).
%$R$ defines connections of boundary points by arcs. 
Consider the connection of $a_0$, that is, the boundary point connected 
to $a_0$ by an arc in $R$. Initially, 
we have two cases: $a_0$ connects
to either $a_m$ or $b_m$ for some $m$.
\par\vs
\begin{center}
\begin{tabular}{cc} 
\includegraphics[trim=0mm 0mm 0mm 0mm, width=.25\linewidth]
{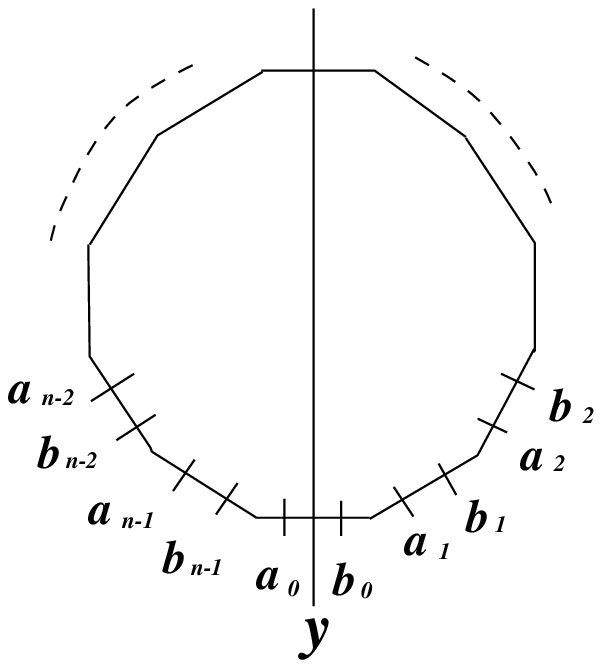}
&\hspace*{20mm}
\includegraphics[trim=0mm 0mm 0mm 0mm, width=.35\linewidth]
{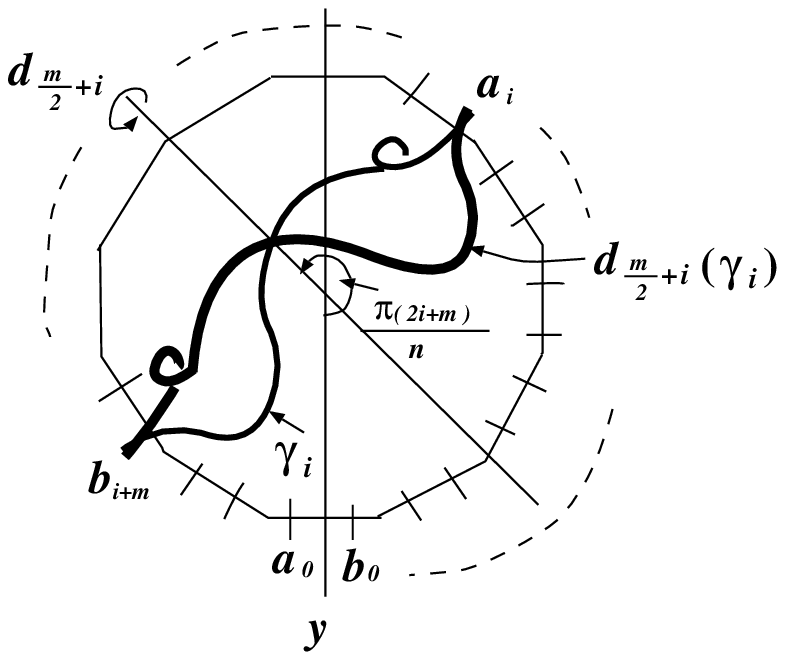}
\\
(a) &\hspace*{20mm} (b)
\end{tabular}
\par
Fig. 2.3
\end{center}

\par
%\begin{enumerate}
%\item[(1)] 
If $n>2$ then $a_0$ cannot be connected to $a_m$. To prove this claim let 
us assume, by  contradiction, 
that $a_0$ connects to $a_m$ then $\phi^m(a_0)=a_m$ connects to
$\phi^m(a_m)=a_{2m}$ which must be the same as $a_0$. Therefore, $2m=n$.
This implies that 
%holds for any boundary points, that is 
$a_i$ connects to $a_{i+ {\frac{n}{2}}}$
and $b_i$ to $b_{i+ {\frac{n}{2}}}$. The arc $\gamma (x_i)$ of $R$ 
connecting $x_i$ to $x_{i+ {\frac{n}{2}}}$, where the symbol $x$ may
	   stand for $a$ or $b$, is setwise preserved
by the rotation ${\phi}^{\frac{n}{2}}$.
% with endpoints exchanged. 
Therefore the arc
	   ${\gamma}(x_i)$ has
one fixed point, namely the point of the intersection with the $z$-axis.
For $n>2$ we have at least two arcs of the type  $\gamma (a_i)$. Such
 arcs cut the $z$-axis at different heights, say $h_i$. On the other hand
$\phi (\gamma (a_i)) = \gamma (a_{i+1})$, so $h_i = h_{i+1}$, which gives
   a contradiction. So in this case, we have $n=2$, and in this case
 Theorem 2.1 follows easily.
\par\vs
%\item[(2)] 
Suppose $a_0$ is connected to $b_m$ for some $m$. Let
	   ${\gamma}_i={\gamma}(a_i)$ denote the arc connecting 
the point $a_i$ with $b_{i+m}$ in $R$. Consider the dihedral 
flype $d_{\frac{m}{2}+i}$ exchanging $a_i$ with $b_{i+m}$. The
	   image $d_{\frac{m}{2}+i}({\gamma}_i)$ connects the same
	   points on the boundary as ${\gamma}_i$ that is $a_i$ and
	   $b_{i+m}$~~(Fig. 2.3(b)), so two boundary points of $R$ are
	   connected in $R$ if and only if they are connected in 
$d_0(R)=d_{{\frac{m}{2}}+i}(R)$. 
In particular, the link $L_1=S \cup R$ and its rotant 
$L_2=S \cup d_0(R)$ have the same number of components.
%\end{enumerate}
\par\vs
By observations similar to the above, we have 
% To show that ${\ell} k(L)$ (for oriented rotants) 
%and $sw(L)$ (for any rotants)
%is preserved it suffices to prove the following claim.
%straightforward but though important Lemma 2.2.  
\par\vs\noindent
%\begin{lem}
%{\bf{Lemma 2.2}}~~~
{\bf{Claim 2.2}}~~~
(i)~{\em {For an unoriented rotor $R$ choose any orientation {\em
 (}directions{\em )} of
its arcs {\em (}e.g. from $a_j$ to $b_{j+m}${\em )}. Let
$I({\gamma}_j,{\gamma}_k)$ denote 
 the sum of sign of crossings ${\gamma}_j$ and ${\gamma}_k$, possibly
 $j=k$, then $I({\gamma}_j,{\gamma}_k) = 
I(d_{\frac{j+k+m}{2}}({\gamma}_j),d_{\frac{j+k+m}{2}}({\gamma}_k))$}}.
\par\noindent
(ii)~{\em{For an oriented rotor $R$ and a closed component $\alpha$ of $R$,
$I({\gamma}_i, {\alpha})=I(d_{{\frac{m}{2}}+i}({\gamma}_i), 
d_{{\frac{m}{2}}+i}({\alpha}))$.}} 
\par\vs
Notice that
$\partial{{\gamma}_j}=\partial(d_{\frac{j+k+m}{2}}({\gamma}_k)),\partial{{\gamma}_k}=\partial(d_{\frac{j+k+m}{2}}({{\gamma}_j}))$ 
and $\partial{{\gamma}_i}=\partial(d_{{\frac{m}{2}}+i}({\gamma}_i)).$
%\end{lem}
\par\vs\noindent
{\bf{Proof}}~~
%\par
(i)~The dihedral flype $d_{\frac{j+k+m}{2}}$ of $R$ sends $a_j$ to $b_{k+m}$ 
and $a_k$ to $b_{j+m}$, thus it sends the arc ${\gamma}_j$, connecting
$a_j$ with $b_{j+m}$ in $R$ (resp. $a_k$ with $b_{k+m}$) to the arc
$d_{\frac{j+k+m}{2}}({\gamma}_k)$ connecting $b_{k+m}$ with $a_k$ in $d_0(R)$
(resp. $d_{\frac{j+k+m}{2}}({\gamma}_j)$ connecting $b_{j+m}$ with $a_k$) 
(Fig. 2.4). Therefore $I({\gamma}_j,{\gamma}_k)=
I(d_{\frac{j+k+m}{2}}({\gamma}_k),d_{\frac{j+k+m}{2}}({\gamma}_j))$, as
required.
\hfill$\Box$ 
\par\vs
\begin{center}
\hspace*{8mm}
\includegraphics[trim=0mm 0mm 0mm 0mm, width=.35\linewidth]
{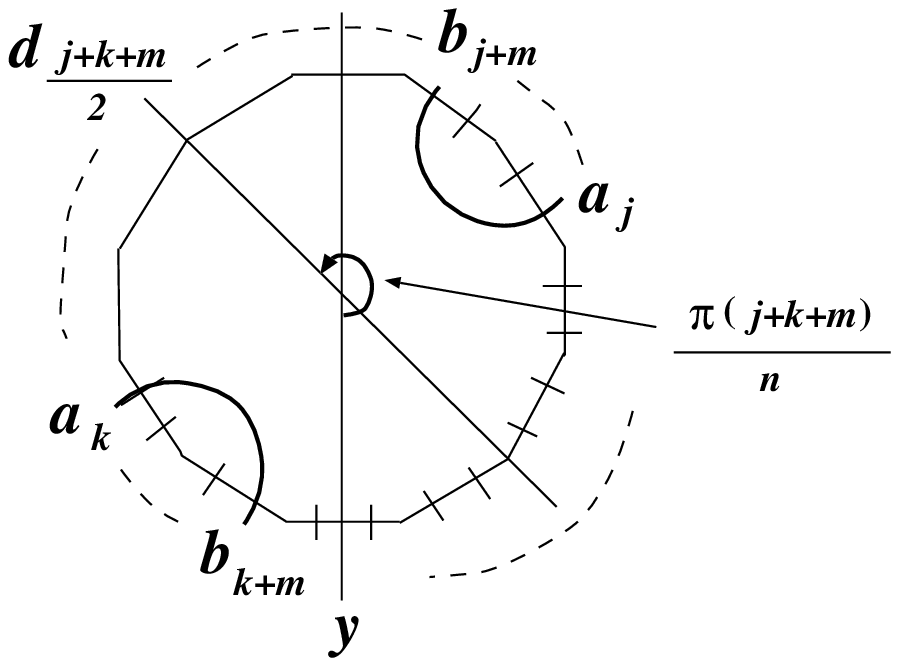}
\\
Fig. 2.4 
\end{center}
\par\noindent
(ii)~Since ${\gamma}_i$ in $R$ and $d_{{\frac{m}{2}}+i}({\gamma}_i)$
in $d_0(R)$ connect the same boundary points $a_i$ and $b_{i+m}$, we
have the conclusion.
\par
Theorem 2.1 (ii), (iii) and (iv) follows from Claim 2.2 and the 
fact that $L_1$ and $L_2$ have the same stator.  
%also 
\hfill$\Box$
\par\vs\vs\vs\vs

We use Theorem 2.1 to show that with some technical assumptions, 
that are explained below, 
the double branched covers of $S^3$ branched along rotant links 
have isomorphic first homology groups. We do not use 
later in the paper the result of  Corollary 2.3, however, we would like 
to contrast it with the example in Fig. 1.1 of rotant links with 
different first homology groups.
% of the double branched covers of $S^3$ branched along them. 

In the proof of Corollary 2.3 we 
use Montesinos method \cite{Mo} 
of finding surgery description of the double branched covers of $S^3$ 
branched along links, when a surface (possibly unoriented) bounding the 
link is given. We closely follow, in this part of the paper, notation 
used in  \cite{IPY}.

Let $T_0$ be a trivial $n$-tangle diagram as in Fig. 2.5(a).
%{\footnote{That is, $T_0$ is obtained from the tangle by a 
%homeomorphism of $D_3$}}.
Let $D_1 \cup \cdots \cup D_n$ be a disjoint union of disks bounded by
$T_0$ and a disjoint union of arcs in $\partial B^3$ connecting $\partial T_0$.
Let $b_1,\cdots,b_m$ be mutually disjoint disks (ribbons) in $B^3$
such that $b_i \cap {\bigcup}_j D_j=\partial b_i \cap T_0$ are two
disjoint arcs in $\partial b_i (i=1,\cdots,m)$, Fig. 2.5(b) . 
We denote by $\Omega (T_0 ; \{D_1,\cdots,D_n\},\{b_1,\cdots,b_m\})$ 
the tangle $T_0 \cup {\bigcup}_i {\partial}b_i -$int$(T_0 \cap 
{\bigcup}_i \partial b_i)$ together with the the surface 
$ \bigcup D_i \cup \bigcup b_j$ and its decomposition into disks 
$D_i$ and $b_i$. 
We call such a structure a {\it{disk-band representation}} 
of a tangle \cite{IPY}. 
\par

\par\vs
\begin{center}
\begin{tabular}{cc} 
\includegraphics[trim=0mm 0mm 0mm 0mm, width=.30\linewidth]
{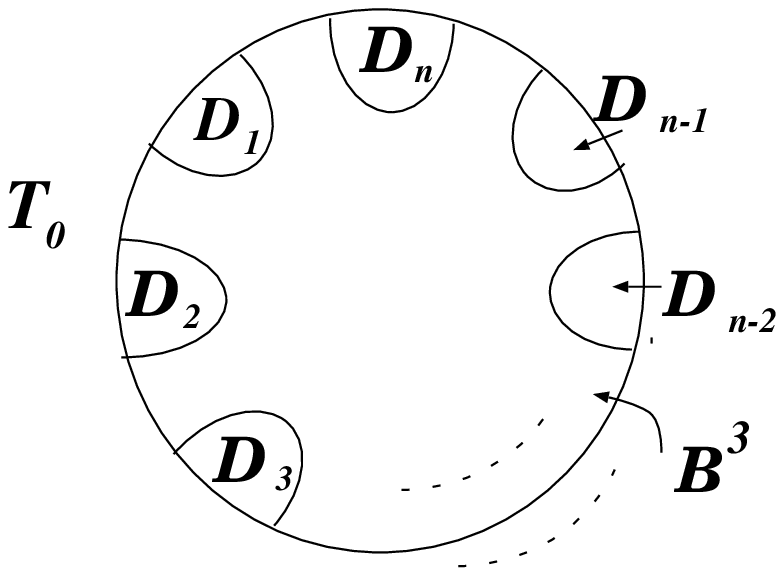}
&\hspace*{20mm}
\includegraphics[trim=0mm 0mm 0mm 0mm, width=.30\linewidth]
{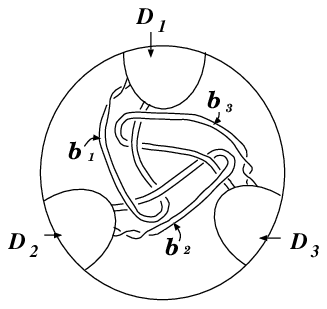}\\
(a)&\hspace*{20mm}$(b)$
\end{tabular}
\\
Fig. 2.5 
\end{center}
\ \\

If a rotor part has a rotationally symmetric disk-band
representation, then the following corollary of Theorem 2.1 holds.
\par\vs\vs\noindent
{\bf{Corollary 2.3}}~~{\it{ Let $L_1$ and $L_2$ be a pair of unoriented 
$n$-rotants such that $n$-rotor $R_1$ of $L_1$ admits a rotational symmetric 
disk-band representation with the number of ribbon disks in the 
representation equal to $n$.  
%Let $L_1$ and $L_2$ be a pair of rotants of any order $n$.
%Assume that $n$-rotors $R_k$ of $L_k$ are in the disk-band representations
%$\Omega(T_0;\{D_1,...,D_n\},\{b_{k1},...,b_{kn}\})$\\
%$\Omega_k(T_0;\{D_1,\cdots,D_n\},\{b_1,\cdots,b_n\})$~~
%$(k=1,2)$ which have also $n$-rotational symmetries. 
%(We assume that  the numbers of the elements of $\{D_1,\cdots,D_n\}$ 
%and $\{b_{k1},\cdots,b_{kn}\}$ are same.)
Then $H_1(M^{(2)}_{L_1},{\Bbb Z})= H_1(M^{(2)}_{L_2};\Bbb Z)$  
where $M^{(2)}_{L}$ 
denotes the double branched cover of $S^3$ branched along a link $L$.
}}
\par\vs
\par\vs\noindent
{\bf Proof} Let $\Omega(T_0;\{D_1,...,D_n\},\{b_{k1},...,b_{kn}\})$ 
be the the disk-band representations of $R_1$ and $R_2=d_0(R_1)$
respectively,
related by the dihedral flype $d_0$.
Let $B^3$ be the $3$-ball such that $B^3\cap L_k$ is the tangle ingredient
 of $\Omega(T_0;\{D_1,...,D_n\},\{b_{k1},...,b_{kn}\})$ 
$(k=1,2)$ and $B^3_0=B^3-{\mathrm int}N(D_1\cup\cdots\cup D_n)$, 
where $N(D_1\cup\cdots\cup D_n)$ is a regular neighborhood 
of $D_1\cup\cdots\cup D_n$ in $B^3$. 
There are compact, connected, possibly non-orientable surfaces $F_k$ $(k=1,2)$ 
in $S^3$ such that $F_k\cap B^3=D_1\cup\cdots\cup D_n\cup 
b_{k1}\cup\cdots\cup b_{kn}$ and the surface $F_k \cap (S^3-B^3)$ is 
connected. We follow \cite{IPY} in constructing a surgery 
description of the 
double branched cover using a surface $F_k$. We work with $F_1$ and and $L_1$,
 the construction for $F_2$ and $L_2$ is related by a dihedral flype.

Choose a point $v_i$ in $D_i\cap \partial B^3$ $(i=1,...,n)$.
Let $G_k$ be a spine of 
$F_k$ with the vertex set $\{v_1,...,v_n\}$ such that 
$G_k\cap B^3$ is a spine of $D_1\cup\cdots\cup D_n\cup 
b_{k1}\cup\cdots\cup b_{kn}$. Let $T_k\subset S^3-{\mathrm int}B^3$
be a spanning tree of $G_k$ and $G_k/T_k$ a spine obtained from 
$G_k$ by contracting $T_k$ into a point $v$. We may assume that 
$N(G_k/T_k)\cap F_k$ consists of a disk $D_{k0}$ containing $v$ 
and mutually disjoint disks $b'_{k1},...,b'_{km}$ such that 
$b'_{ki}\cap D_{k0}= \partial b'_{ki}\cap\partial D_{k0}$ are 
two disjoint arcs in $\partial b'_{ki}$ $(i=1,...,m)$, 
$(D_{k0}\cup b'_{k1}\cup\cdots\cup b'_{km})\cap B^3_0=
(b_{k1}\cup\cdots\cup b_{kn})\cap B^3_0$, and that 
$(D_{10}\cup b'_{11}\cup\cdots\cup b'_{1m})-B^3_0=
(D_{20}\cup b'_{21}\cup\cdots\cup b'_{2m})-B^3_0$. 
Let $\varphi:S^3\rightarrow S^3$ be the double branched cover 
branched along $\partial D_{k0}$. 
Then $M^{(2)}_{L_k}$ is obtained from $S^3$ by surgery along a 
framed link $\varphi^{-1}(b'_{k1}\cup\cdots\cup b'_{km})$. 
Note that $\varphi^{-1}((b'_{k1}\cup\cdots\cup b'_{km})\cap B^3_0)=
\varphi^{-1}((b_{k1}\cup\cdots\cup b_{kn})\cap B^3_0)$ are two $n$-rotors 
and $\varphi^{-1}((b'_{11}\cup\cdots\cup b'_{1m})- B^3_0)=
\varphi^{-1}((b'_{21}\cup\cdots\cup b'_{2m})- B^3_0)$. 
Since each $\varphi^{-1}(b'_{ki})$ is a component of 
$\varphi^{-1}(b'_{k1}\cup\cdots\cup b'_{km})$,  
it is not hard to see that there is a blackboard framed, oriented link
$c_{k1} \cup \cdots \cup c_{km}$ such that each $c_{ki}$ corresponds to
$b'_{ki}$ and the both components of $(c_{k1} \cup \cdots \cup c_{km})
\cap \varphi^{-1}(B^3_0)$ are oriented $n$-rotors.
So $c_{21} \cup \cdots \cup c_{2m}$ is obtained from 
$c_{11}\cup\cdots\cup c_{1m}$ by 
%twice 
two oriented 
$n$-rotations. By Theorem 2.1 (iii), the linking matrices of 
$c_{11}\cup\cdots\cup c_{1m}$ and 
$c_{21}\cup\cdots\cup c_{2m}$ coincide. 
Since the linking matrix of $c_{k1}\cup\cdots\cup c_{km}$
is a relation matrix of the first homology group of $M^{(2)}_{L_k}$, 
we have the conclusion. 

\hfill$\Box$

Corollary 2.3 and the example in Fig 1.1 allow us to conclude that 
not every $n$-rotor has a symmetric disk-band representation with 
$n$ bands. 

%\hfill$\Box$
%\begin{df}
\par\vs
\par\vs\vs
Let $F_L$ be a Seifert surface of an oriented link $L$.
Denote by $\psi : H_1(F_L;\Bbb Z$$) \times H_1(F_L;\Bbb Z$$) \rightarrow
\Bbb Z$ 
the Seifert form associated with $F_L$ (i.e. ${\psi}(x,y)=lk(x^+,y)$,
 where $x^+$ denotes a curve pushed
$x$ slightly off $F_L$ into the positive direction).
Choosing an ordered
basis for $H_1(F_L;\Bbb Z)$ allows us to describe 
the form $\psi$ by the corresponding
 Seifert matrix. Let ${\cal A}_L$ be the Seifert 
matrix of the form $\psi$ with respect to some ordered basis of
$H_1(F_L;\Bbb Z)$. 
\par\vs\vs
Let $F_L$ be a spanning surface, possibly nonorientable, of an
unoriented link $L$. We use the following generalization of 
Seifert\footnote{It is a generalization of the symmetrization of 
the Seifert form.}  and
Goeritz forms defined by 
Gordon and Litherland in \cite{G-L}. For the spanning surface $F_L$
 consider regular neighborhood, $N(F_L)$, of $F_L^3$ in $S^3-L$. 
Then $N(F_L)$ is 
the $I$-bundle over $F_L$ and the $\partial I$-bundle $\tilde{F_L}$
is a double cover of $F_L$ (possibly disconnected) with 
the projection map $p: \tilde{F_L}
\rightarrow F_L$.  The bilinear form ${\cal G}_{F_L}: H_1(F_L;\Bbb Z$$)
\times H_1(F_L;\Bbb Z$$) \rightarrow \Bbb Z$ defined by 
${\cal G}_{F_L}(x,y)=lk(p^{-1}x,y)$, where $x$ and $y$ are oriented
loops in $F_L$, is called the {\it{Goeritz form}} associated to the surface
 $F_L$. 
%Any choice of 
For an ordered basis of $H_1(F_L;\Bbb Z)$ we have
%determines 
the matrix  $G_{F_L}$ representing the Goeritz form ${\cal G}_{F_L}$. 
The matrix $G_{F_L}$ is called the {\it{Goeritz matrix}} of $F_L$ 
with respect to a basis of $H_1(F_L;\Bbb Z)$.
% associated with $F_L$.
%\par\vs\vs
% and $G_{F_L}$ a
% Goeritz matrix.
%Let $L'=K'_1 \cup K'_2 \cup
% \cdots \cup K'_m$ be a $m$-component link  
%which is parallel to $L$
% on the surface $F_L$. Suppose 
%the orientations of each
% pair $K_i$ and $K'_i$ are parallel.
% Then we define the normal Euler
% number $e(F_L)=- {\displaystyle{\sum_{i=1}^{m}}} lk (K_i,K'_i)$.
\par
 The form ${\cal G}_{F_L}$ defined 
over $\Bbb Z$ can be extended to the form ${\hat{\cal G}}_{F_L}$ over $\Bbb C$. We
view the form ${\hat{\cal G}}_{F_L}$ as the Hermitian form represented in a  
%any 
basis by
the Hermitian matrix ${\hat{G}}_{F_L}$~~(i.e. 
${\hat{G}}_{F_L}={\overline{{\hat{G}}^t_{F_L}}}$).
\par
\vs\vs
For a spanning surface $F_{L_k}$ of $L_k=K_{k1} \cup K_{k2} 
\cup \cdots \cup K_{km}$, the framing of $L_k$ is uniquely 
determined by $F_{L_k}$ as follows\footnote{The regular neighborhood 
of $K_{ki}$ in $F_k$ is the frame knot associated to $K_{ki}$. 
Its framing, when compared to the standard framing, is given by 
$lk(K_{ki},{K_{ki}}^{F_{L_k}})$.}: 
Let ${K_{ki}}^{F_{L_k}}$ be a parallel copy of $K_{ki}$ that misses 
$F_{L_k}$. We define the framing 
$K_{ki}$ to be $lk(K_{ki},{K_{ki}}^{F_{L_k}})$. 
We put $e(F_{L_k})=-{\sum}_{i}^{} lk (K_{ki},{K_{ki}}^{F_{L_k}})=- tr(L_k)$. 
%The following corollary is a consequence of Theorem 2.1.
%\par\vs\noindent
%\begin{cor}
%{\bf{Corollary 2.7}}~~~{\em{If $L_2$ is obtained from $L_1$ 
%by a rotation, then $ e(F_{L_1})=e(F_{L_2})$.}}
%\hfill$\Box$
%\end{cor}
\par\vs\vs
We recall the definition of the Tristram-Levine signature of an oriented link.
\par\vs\noindent
%\begin{df}
{\bf{Definition 2.4}}~~\cite{L,Tr}~~~{\em{Let $L$ be an oriented link in
 $S^3$ and let $\omega$ be a complex number with $\mid \omega
 \mid =1$, $\omega \ne 1$. The Tristram-Levine signature of $L$, denoted
 by ${\sigma}_{\omega}(L)$, is the signature of the
 Hermitian matrix $(1-\omega){\cal A}_L+(1-\bar{\omega}){\cal A}_L^t$,
 where ${\cal A}_L$ is
 a Seifert matrix of $L$. }}
%\end{df}
\par\vs\noindent
%\begin{df} 
{\bf{Definition 2.5}}~~\cite{G-L,M,M2}~ {\em{
Let $L$ be an unoriented link in $S^3$, and let ${\hat{L}}$ be the link
 obtained from $L$ by a choice of an orientation.
The Murasugi signature ${\hat{\sigma}}(L)$ of an unoriented link $L$ 
is defined to be ${\hat{\sigma}}(L)=\sigma({\hat{L}})+ {\ell}k({\hat{L}})$.}}
%\end{df}
\par\vs\noindent
%\begin{remark}
{\bf{Remark 2.6}}~~~
Murasugi showed  in \cite{M2} that $\sigma (\hat{L})+ 
 {\ell}k({\hat{L}})$ does not depend on the choice of an orientation 
of $L$. So $\hat{\sigma}(L)$ is an invariant of unoriented links. 
We shall use later the fact that
 $\hat{\sigma}(L)=$sign$(G_{F_L}) + {\frac{1}{2}}e(F_L)$ \cite{G-L}. 
%\end{remark}
\par\vs

\par\vs
\section{Unoriented rotation and Murasugi signature}
%\section{Rotation and double cover of links}
~~~\par\vs
In this section we prove that the Murasugi signature of unoriented links 
%\cite{G-L, M, M2} 
is preserved by any rotation. 
The result follows from a more general statement (see Theorem 3.2) 
%We place this result in a more general setting by proving 
that the rotation preserves the characteristic polynomial
of the Goeritz matrix (with the special choices of surfaces). 
In particular Theorem 3.2 allows us to obtain the result mentioned 
first in \cite{P} that was also proven by Traczyk that the determinant of an
unoriented link is preserved by any rotation.
\par\vs
\begin{theorem}
Let $L_1$ and $L_2$ be a pair of unoriented $n$-rotants
 (no restrictions on $n$). 
Then $\hat{\sigma}(L_1)=\hat{\sigma}(L_2)$.
\end{theorem}
\par\vs
The main result of this section is Theorem 3.2 from which Theorem 3.1 follows.
\par
\par\vs
Let $L_1$ and $L_2$ be a pair of unoriented rotant links.
Consider projections of the links $L_1$ and $L_2$ onto $\Bbb R$$^2$ with
rotor parts $R_1$ and $R_2$ contained in disks $D_1$ and $D_2$, 
respectively. We
can deform the stator parts $S_1$ and $S_2$ of the diagrams of $L_1$ and 
$L_2$ into the position shown in Figure 3.1.
\par
We color the regions on $\Bbb R$$^2$ bounded by the diagrams of $L_k$
in a checkerboard manner as in Figure 3.2. Using the
black regions we form the spanning surface $F_{L_k}$ for $k=1,2$. 
We choose for a basis of $H_1(F_{L_k};\Bbb Z)$ the anti-clockwise
oriented boundary curves of the bounded white regions, and we refer 
to this basis as the {\it{standard basis}}.  
\par\vs
\begin{center}
\begin{tabular}{cc} 
\includegraphics[trim=0mm 0mm 0mm 0mm, width=.5\linewidth]
{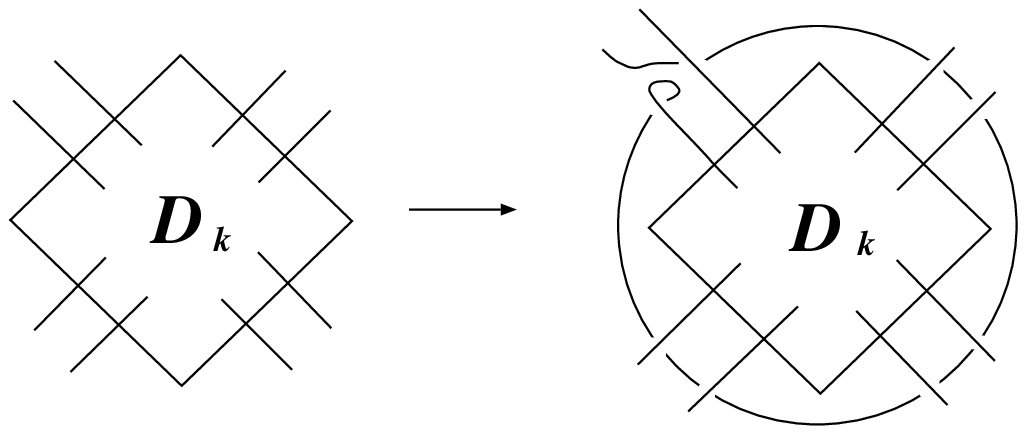}
&\hspace*{20mm}
\includegraphics[trim=0mm 0mm 0mm 0mm, width=.25\linewidth]
{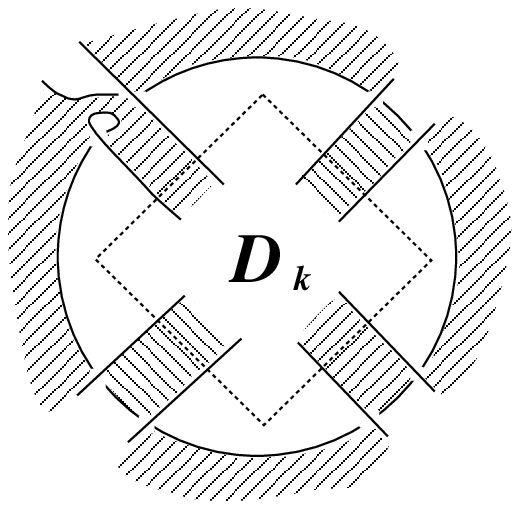}
\\
Fig. 3.1 &\hspace*{20mm} Fig. 3.2
\end{tabular}
\end{center}
\par
We also may assume that the framed links $L_1$ and $L_2$ obtained 
from $F_1$ and $F_2$ respectively, form a pair of
rotants.
By Theorem 2.1, $tr(L_1)=tr(L_2)$, so we have $e(F_{L_1})=e(F_{L_2})$.
This fact, Remark 2.6 and the following theorem imply Theorem 3.1.
\par\vs
With the choices for $F_k$'s and bases of $H_1(F_k;\Bbb Z)$'s, 
made above, we can formulate the main result of this section.
\par\vs
\noindent
{\bf{Theorem 3.2}}~~
%\begin{theorem}~~
{\it{Let $G_{F_k}~(k=1,2)$ be the Goeritz matrices with respect to the
standard basis. Then $\det(G_{F_{L_1}}-{\lambda}E)= \det(G_{F_{L_2}}-{\lambda}E)$}}.
%The characteristic polynomials of Goeritz matrices of
%a pair of rotants coincide, that is 
%\end{theorem}
\par\vs
%\noindent
{\bf{Proof}}\footnote{We adjust here the Traczyk's method \cite{T} to the 
case of unoriented rotors and Goeritz matrices.}~~
%\par
Let $X_{{\cal S}_k}$ and $X_{{\cal R}_k}$ be the subsets
of the standard basis of $H_1(F_{L_k};\Bbb Z)$ which live entirely in 
the stator and rotor part respectively, and let $X_{{\cal M}_k}$ be the
complement of $X_{{\cal S}_k} \cup X_{{\cal R}_k}$ in the standard basis.
$X_{{\cal M}_k}$ is composed of boundaries of white regions intersecting 
the boundary of the rotor. We can have $n$ such regions or just one region.
We can, however, always assume, modifying the rotor part of
the diagram if necessary, that $X_{{\cal M}_k}$ has $n$ different elements.
Consider submodules ${\cal S}_k$, ${\cal R}_k$ and ${\cal
M}_k $ of $H_1(F_{L_k};\Bbb Z)$ generated by $X_{{\cal S}_k}$, $X_{{\cal R}_k}$ and $X_{{\cal M}_k}$.
We have the following
decomposition into the direct sum of $\Bbb Z$-modules : 
$H_1(F_{L_k};\Bbb Z)$$={\cal S}_k \oplus {\cal M}_k \oplus {\cal R}_k$. 
%The common part of ${\cal S}_k$ and ${\cal M}_k+{\cal R}_k$ 
%is generated by the ``circle ''$C_k$ (Fig.3.2). 
Let $v$ denote the generator of ${\cal M}_1$ intersecting the $y$ axis 
of the dihedral flype $d$ (Fig. 3.3). 
% We can always assume, modifying the rotor part of 
%the diagram if necessary, that the boundaries 
%of middle regions $v, {\alpha}(v),\cdots,{\alpha}^{n-1}(v)$
% form an ordered basis of ${\cal M}$.
%Deletion of one standard generator from $X_{{\cal S}_k}$ 
%gives the new set $X_{{\cal S}'_k}$ which generates the submodule 
%$S'_k$ of $H_1(F_{L_k};\Bbb Z)$. 
%\par
%We have the split of the homology into 
%the direct sum $H_1(F_k;\Bbb Z)={\cal S}'_k \oplus {\cal
%M}_k \oplus {\cal R}_k$ where $X_{{\cal S}'_k}$, $X_{{\cal
%R}_k}$ and $X_{{\cal M}_k}$ are the bases of ${\cal S}'_k$,
%${\cal R}_k$ and ${\cal M}_k$ respectively. 
There is an action of the cyclic group 
$\Bbb Z$$_n=< \alpha\ |\ \alpha^n=1 >$ on ${\cal R}_1 \oplus {\cal M}_1$ 
induced by the
$\frac{2\pi}{n}$-rotation around the center of $D_1$. Thus the ordered set
$X_{{\cal M}_1}: v, 
\alpha(v),{\alpha}^2(v),\cdots,{\alpha}^{n-1}(v)$ can be assumed to be  
%a generating set
a basis
%\footnote{We can assume without lost of generality that it is a basis.
%Otherwise the middle set has only one white 
%region and the theorem is obvious in this case as rotants differ only by 
%attachment of connected sums.}
of ${\cal M}_1$.
Let $X^{\ast}_{{\cal R}_1}$ be a set of generators of ${\cal R}_1$
formed by choosing one representative from each 
orbit of $\Bbb Z$$_n$-action on standard
generators of ${\cal R}_1$ (i.e. $X^{\ast}_{{\cal R}_1}=X_{{\cal R}_1}/
\Bbb Z$$_n$). We construct a bijection ${\eta}$ between the set of standard
generators of $H_1(F_{L_1};\Bbb Z)$ and $H_1(F_{L_2};\Bbb Z)$. First, define 
$\eta {\mid}_{X_{{\cal S}_1}} : X_{{\cal S}_1} \rightarrow X_{{\cal S}_2}$ 
to be the identity map since the stator part is unchanged by rotation. 
The map ${\eta} {\mid}_{X_{{\cal M}_1}} : X_{{\cal M}_1}
\rightarrow X_{{\cal M}_2}$ is given by ${\eta}
({\alpha}^j(v))={\alpha}^j(d(v))$~~(i.e. ${\alpha}^j(v)$ and
${\eta}({\alpha}^j(v))$ have the same stator parts). 
Finally, ${\eta} {\mid}_{X_{{\cal R}_1}}: X_{{\cal R}_1} \rightarrow X_{{\cal 
R}_2}$ is given by  ${\eta}({\alpha}^j(x))=d({\alpha}^j(x))$ for $x \in
X_{{\cal R}_1}$. The bijection ${\eta}$ extends to the isomorphism, 
$H_1(F_{L_1};\Bbb Z)$ $\rightarrow H_1(F_{L_2};\Bbb Z)$, that is 
also denoted by ${\eta}$.
We use the isomorphism $\eta$ to identify $H_1(F_{L_1};\Bbb Z)$ with
$H_1(F_{L_2};\Bbb Z)$. This identification allows us to drop the indices 
in ${{\cal S}}_k, {\cal M}_k$ and ${\cal R}_k$ 
and write ${\cal S},{\cal M}$ and ${\cal R}$. 
\par\vs
\begin{center}
\begin{tabular}{cc} 
\includegraphics[trim=0mm 0mm 0mm 0mm, width=.35\linewidth]
{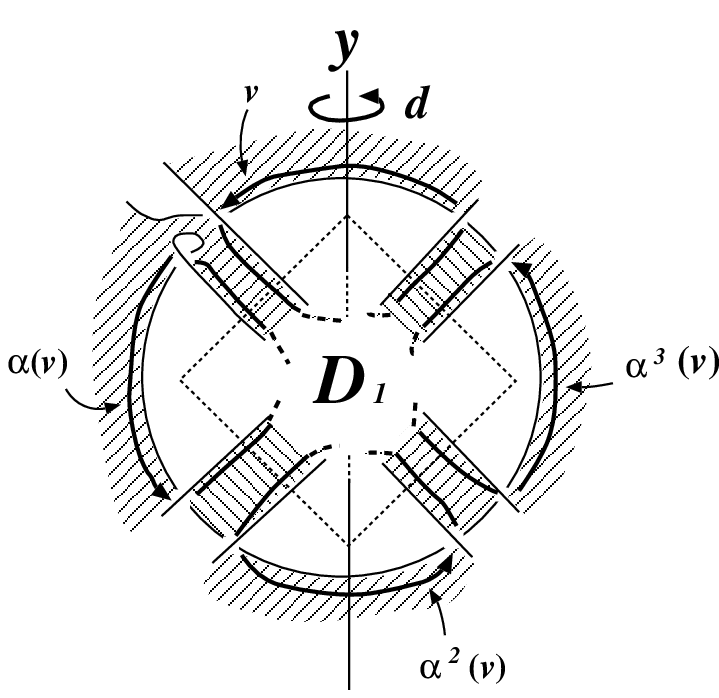}
&\hspace*{20mm}
\includegraphics[trim=0mm 0mm 0mm 0mm, width=.35\linewidth]
{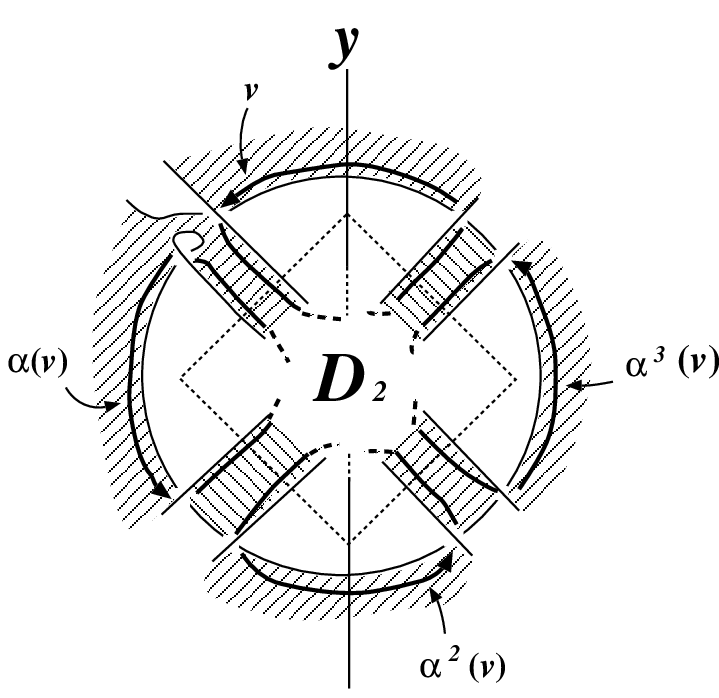}
\end{tabular}
\\
Fig.3.3
\end{center}
\par\vs
Let us consider forms 
${\cal G}_1={\cal G}_{F_{L_1}}$ and ${\cal G}_2={\cal G}_{F_{L_2}}$ on
the same space ${\cal S} \oplus {\cal M} \oplus {\cal R}$.
\par
We have the following properties of ${\cal G}_1$ and ${\cal G}_2$.
%(0).~~${\hat{{\cal G}}}_k(p,q)={\hat{{\cal G}}}_k(q,p)$ 
%for all $p,q \in {\cal S}\oplus{\cal M}+{\cal R}$
\par 
\begin{enumerate}
\item[(1)]  ${\cal G}_2 (x,y) ={\cal G}_1(x,y)$ for all $ x,y \in {\cal S} \oplus {\cal M}$.
\item[(2)]  ${\cal G}_2 (x,y) ={\cal G}_1(x,y)$ for all
generators $ x,y \in {\cal R}$ and
%\item[~~]  ${\hat{\cal {G}}}_2 (x,y) = {\hat{\cal G}}_1(x,y)=
%{\overline{{\hat{\cal G}}_1({\bar{x}},{\bar{y}})}}=
%{\hat{\cal G}}_1({\bar{y}},{\bar{x}})$ for all $ x,y \in {\bf{R}}$.
\item[(3)]  ${\cal G}_1(x,y)={\cal G}_1({\alpha}^l(x),{\alpha}^l(y))$ for all 
generators $x,y \in {\cal M} \oplus {\cal R}$, 
\item[~~]  ${\cal G}_2(x,{\alpha}^l(v))={\cal G}_1(x,{\alpha}^{-l}(v))$ for every 
generator $x$ of ${\cal R}$,
\item[~~]  ${\cal G}_2({\alpha}^l(x),v)={\cal G}_1({\alpha}^l(x),v)$ for every 
generator $x$ of ${\cal R}$, and
\item[~~]  ${\cal G}_2(x,v)={\cal G}_2({\alpha}^l(x),{\alpha}^{-l}(v))$ for every 
generator $x \in {\cal R}$.
\item[(4)]  ${\cal G}_k(x,y)=0$ for all $x \in {\cal S}, 
y \in {\cal R}, (k=1,2)$.
%\item[(5)]  ${\cal G}_k(x,y)={\cal G}_k(y,x)$ 
%for all $x,y  \in {\cal S}' \oplus {\cal M} \oplus {\cal R},(k=1,2)$.
\end{enumerate}
\par\vs
Let ${\bf{S}},{\bf{M}}$ and
${\bf{R}}$ be the subspaces of $ ({\cal S}\oplus {\cal M}\oplus {\cal R})
\otimes \Bbb C$ complexifying 
${\cal S}, {\cal M}$ and ${\cal R}$, respectively. We have the 
involution ~${\bar{}}$~ : ${\bf{S}}\oplus {\bf{M}} \oplus {\bf{R}} 
\rightarrow {\bf{S}}\oplus {\bf{M}} \oplus {\bf{R}} $ corresponding to 
the conjugation in the factor $\Bbb C$ of the
tensor product. The image of $x \in {\bf{S}} \oplus {\bf{M}} \oplus
{\bf{R}}$ 
under this
involution is denoted by ${\bar{x}}$.
%For the Hermitian matrices ${\hat{G}}_{L_k}~~(k=1,2)$,
%we may choose a basis of $H_1(F_{L_k};C)$, from generators of 
%$H_1(F_{L_k};Z)$. 
Using the rotational symmetry of the rotor part we conveniently change 
the basis of ${\bf{M}}$ and the generating set of ${\bf{R}}$ 
in the following way. 
Let ${\omega}_j$ be an $n$th root of unity, 
${\omega}_j=e^{2 \pi i {\frac{j}{n}}}$.
We replace the basis $\{{\alpha}^j(v) \mid j=0,1,\cdots, n-1\}$
of ${\bf{M}}$ by $\{{\bf{v}}_j \mid
 {\bf{v}}_j={\displaystyle{\sum_{l=0}^{n-1}}}{\omega}_j^l{\alpha}^l(v),
 j=0,1,\cdots,n-1\}$. For ${\bf{R}}$ we consider two choices of generating 
sets that are related by the involution 
~${\bar{}}$~ as follows. We either replace  the set $ \{{\alpha}^j(y_p)
 \mid y_p \in X^{\ast}_{\cal R}, j=0,1,\cdots,n-1\}$, by $\{{\bf{y}}_{j,p} \mid
{\bf{y}}_{j,p}={\displaystyle{\sum_{l=0}^{n-1}}}{\omega}_j^l{\alpha}^l(y_p), 
y_p \in X^{\ast}_{{\cal R}},
j=0,1,\cdots,n-1\}$ or by $\{{\overline{{\bf{y}}_{j,p}}} \mid
{\overline{{\bf{y}}_{j,p}}}={\displaystyle{\sum_{l=0}^{n-1}}}{\overline{{\omega}_j}}^{l}{\alpha}^l(y_p), 
y_p \in X^{\ast}_{{\cal R}},j=0,1,\cdots,n-1\}$.
\par
Let us consider the Hermitian forms
${\hat{{\cal G}}}_1={\hat{{\cal G}}}_{F_{L_1}}$ and 
${\hat{{\cal G}}}_2={\hat{{\cal G}}}_{F_{L_2}}$, induced by ${{\cal G}}_1$ 
and ${\cal G}_2$, on
the same space ${\bf S} \oplus {\bf M} \oplus {\bf R}$.
%\par
%We have the following relations for 
%${\hat{{\cal G}}}_1$ and ${\hat{{\cal G}}}_2$.
\par\vs
These new generating sets for ${\bf{M}} \oplus {\bf{R}}$ satisfy 
the following conditions.
\par\medskip 
\begin{enumerate}
\item[(1)]  ${\hat{{\cal G}}}_k({\bf{v}}_j,{\bf{v}}_m)=0$ 
for $j \ne m$, where $ {\bf{v}}_j,{\bf{v}}_m  \in {\bf{M}} $ and
	    $k=1,2$,
\item[~~]  ${\hat{{\cal G}}}_1({\bf{x}}_{j,p},{\bf{v}}_m)
={\hat{{\cal G}}}_2({\overline{{\bf{x}}_{j,p}}},{\bf{v}}_m)=0$ 
for $j \ne m$ where $ {\bf{x}}_{j,p} \in {\bf{R}}_1,{\overline{{\bf{x}}_{j,p}}} \in {\bf{R}}_2, {\bf{v}}_m \in {\bf{M}} $,
\item[~~]  ${\hat{{\cal G}}}_1({\bf{x}}_{j,p},{\bf{y}}_{m,q})
={\hat{{\cal G}}}_2({\overline{{\bf{x}}_{j,p}}},{\overline{{\bf{y}}_{m,q}}})$ 
for $j \ne m$ where $ {\bf{x}}_{j,p},{\bf{y}}_{m,q} \in {\bf{R}}_1,{\overline{{\bf{x}}_{j,p}}},{\overline{{\bf{y}}_{m,q}}} \in {\bf{R}}_2.$
\item[(2)]  ${\hat{{\cal
G}}}_1({\bf{x}},{\bf{y}}_{j,p})={\hat{{\cal
G}}}_2({\bf{x}},{\overline{{\bf{y}}_{j,p}}})=0~~$ 
for any $ {\bf{x}} \in {\bf{S}}, {\bf{y}}_{j,p} \in
	    {\bf{R}}_1,{\overline{{\bf{y}}_{j,p}}} \in {\bf{R}}_2$.
\item[(3)]  ${\hat{{\cal G}}}_1({\bf{y}}_{j,p},{\bf{y}}_{j,q})=
{\overline{{\hat{{\cal G}}}_2({\overline{{\bf{y}}_{j,p}}},{\overline{{\bf{y}}_{j,q}}})}}$, for any ${\bf{y}}_{j,p},{\bf{y}}_{j,q} \in {\bf{R}}_1, {\overline{{\bf{y}}_{j,p}}}, {\overline{{\bf{y}}_{j,q}}} \in {\bf{R}}_2$.
\item[(4)]  ${\hat{{\cal G}}}_1({\bf{v}}_j,{\bf{y}}_{j,p})
={\overline{{\hat{{\cal G}}}_2({\bf{v}}_j, {\overline{{\bf{y}}_{j,p}}})}}$ 
 for any ${\bf{v}}_j \in {\bf{M}}, {\bf{y}}_{j,p} \in  {\bf{R}}_1, {\overline{{\bf{y}}_{j,p}}} \in {\bf{R}}_2$.
\item[~~]  ${\hat{{\cal G}}}_1({\bf{v}}_j,{\bf{v}}_j)
={\hat{{\cal G}}}_2({\bf{v}}_j,{\bf{v}}_j)$
 for any ${\bf{v}}_j \in {\bf{M}}$.
\item[(5)]  ${\hat{\cal G}}_1({\bf{x}},{\bf{v}}_j)=
{\hat{\cal G}}_2({\bf{x}},{\bf{v}}_j)~~$ for any ${\bf{x}} \in {\bf{S}}, {\bf{v}}_j
\in {\bf{M}}$.
\end{enumerate}
\par
\medskip
For a given ${\omega}_j$, $0 \leq j \leq n-1$, let $W_j$ be the 
subspace of ${\bf{M}} \oplus {\bf{R}}$ defined by choosing its 
ordered basis in the following way. Take 
${\bf{v}}_j$ from ${\bf{M}}$ first and ${\bf{y}}_{j,p}$ from ${\bf{R}}$
in any order. To obtain the ordered basis of ${\bf{M}} \oplus
{\bf{R}}$ we place the basis of $W_j$ before the basis of $W_{j+1}$ for
$j=0,1, \cdots, n-1$. Finally we add the ordered basis of ${\bf{S}}$. 
We obtain, in this way, an ordered basis of $H_1(F_{L_1};\Bbb C)$.
Notice that we can construct an ordered basis of 
$H_1(F_{L_2};\Bbb C)$ by replacing each
${\bf{y}}_{j,p}$ with ${\overline{{\bf{y}}_{j,p}}}$. 
\par
Let ${\hat{G}_1}$ and ${\hat{G}_2}$ be the matrices 
of the forms ${\hat{{\cal G}}_1}$  and ${\hat{{\cal G}}_2}$ respectively,
 in the ordered bases of ${\bf{S}} \oplus {\bf{M}} \oplus {\bf{R}}$, chosen 
before. 
\par
\begin{center}
$ {\hat{G}}_{L_1} = \left(
        \begin{array}{cccc}
        B_{10} &  &  {\bf{0}} &{}^t{\overline{S_0}} \\
               &\ddots&  &  \vdots \\
        {\bf{0}} &  & B_{1{n-1}}&   {}^t{\overline{S_{n-1}}} \\
        S_0& \cdots & S_{n-1} & S
   \end{array}
       \right)$ and 
$ {\hat{G}}_{L_2} = \left(
        \begin{array}{cccc}
        B_{20} & & {\bf{0}} & {}^t{\overline{S_0}} \\
               &\ddots& &  \vdots\\
       {\bf{0}} & & B_{2{n-1}} & {}^t{\overline{S_{n-1}}} \\
         S_0 & \cdots  & S_{n-1} & S
   \end{array}
       \right)$.
\par\vspace{0.3cm}
\end{center}
In these bases, $B_{1j}$ (respectively $B_{2,j}$), where $j=0,1,\cdots,n-1$,
is the matrix of the restriction of the form ${\hat{\cal G}}_1$ (and 
${\hat{\cal G}}_2$ respectively) to the subspace $W_j$ generated by 
$\{{\bf{v}}_j\} \cup \{{\bf{y}}_{j,p} \mid y_p \in X^{\ast}_{{\cal
R}_1}\}$ ($\{{\bf{v}}_j \} \cup \{ {\overline{{\bf{y}}_{j,p}}} \mid y_p \in
X^{\ast}_{{\cal R}_1}\}$ respectively). 
Finally, the restrictions of ${\hat{{\cal G}_1}}$ and ${\hat{{\cal G}_2}}$ 
to the stator part
 are the same for ${\hat{{\cal G}_1}}$ and ${\hat{{\cal G}_2}}$ and 
denoted by $S$. 
Notice that ${B_{1k}}^t=B_{2k}$, $S_l=({\bf{s}}_{l1}$~${\bf{0}}$~$\cdots$~${\bf{0}})$, and 
 ${\bf{s}}_{l1}$ is the first column of each matrix $S_l$. 
\par
The matrices $M_k=({\hat{G}}_{L_k} - {\lambda}E)~~(k=1,2)$ satisfy the
conditions of Traczyk's Proposition 2.9 for any real number ${\lambda}$,
\cite{T}.
Thus $\det(M_1)=\det(M_2)$ for any real $\lambda$. So the determinants
are equal for any complex
$\lambda$ as well.
%{\footnote{Eigenvectors 
%of ${\hat{G}}_{L_k}$ restricted to $M$ and
%${{\hat{G}}_{L_k}}-{\lambda}E$ restricted to $M$ are the same
%but they correspond to different eigenvalues : ${\omega}$ in the
%${{\hat{G}}_{L_k}}$ case, and ${\omega}-{\lambda}$ in the 
%${\hat{G}}_{L_k}-{\lambda}E$ case. }}
\hfill$\Box$

\par\vs
%\section{Proof of Theorem 1.2}
\section{Oriented rotation and Tristram-Levine signature}
\par\vs
~~\par
In this section we extend the method developed by Traczyk in {\cite{T}} 
in order to show that orientation-preserving rotations (see Fig. 2.2(a)) 
preserve Conway polynomial. We show that the characteristic 
polynomial of the Hermitian form associated
with the Seifert form of appropriately chosen Seifert surface is 
invariant under orientation-preserving rotations. In particular we 
prove the following result.
\par\vs\noindent
\begin{theorem}
 ~Let $L_1$ and $L_2$ be a pair of orientation-preserving $n$-rotants. 
Then ~$\sigma_{\omega}(L_1)=\sigma_{\omega}(L_2)$.
\end{theorem}
\par\vs
\par
The main result of this section is Theorem 4.2 from which 
Theorem 4.1 follows.  
\par\vs\vs
Let $S^2$ be the sphere of a projection of a link $L$, 
and $F_L$ the Seifert surface of $L$ obtained from the diagram of $L$  
by the Seifert algorithm. Let $H$ be a trivalent 
graph that consist of the Seifert circles and the cores of the
bands. Let $R_1,R_2,\cdots,R_m$ be the components of $S^2-H$ which are
not bounded by Seifert circles. 
Assign the anti-clockwise orientation to each boundary curve of 
the regions $R_i (i=1,\cdots,m)$; then
these curves are generators of $H_1(F_L;\Bbb Z)$. Whenever
we refer to generators of $H_1(F_L;\Bbb Z)$, we mean this particular set of
standard generators for Seifert surface $F_L$.
\par\vs
Let $L_1$ and $L_2$ be a pair of orientation-preserving $n$-rotant 
diagrams.
\par
\medskip
We deform the diagrams $L_k~~(k=1,2)$ on $S^2$ into the position
for which our computation is
feasible, as it was done in \cite{T}. 
Let $D_k$ be a disk in $S^2$ such that $D_k^r=D_k \cap L_k$ 
is the rotor part of the diagram $L_k~(k=1,2)$, 
and ${D^s_k}=\bar{D_k} \cap L_k$ the stator part 
($\bar{D_k}=S^2 - $int$ D_k)$. 
Rotors and stators constructed above are all $n$-tangles. 
We deform the stator part $D_1^s=D_2^s$ to the form shown in Fig. 4.1. 
By doing so we obtain an outermost
Seifert circle $C$ in ${\bar{D}}$ that is parallel to $\partial
{D_k}$. Let ${\overline{D_C}}$ be the region which is bounded by $C$ and
$\partial D_k$ in ${\bar{D}}$.
We extend the rotational symmetries of
the rotor parts $D^r_k~~(k=1,2)$ to the parts embedded in
$D_k \cup {\overline{D_C}}$, i.e., we may assume that 
$D_k \cup {\overline{D_C}}\ (k=1,2)$ contain $n$-rotors.  
\begin{center}
\begin{tabular}{cc} 
\includegraphics[trim=0mm 0mm 0mm 0mm, width=.5\linewidth]
{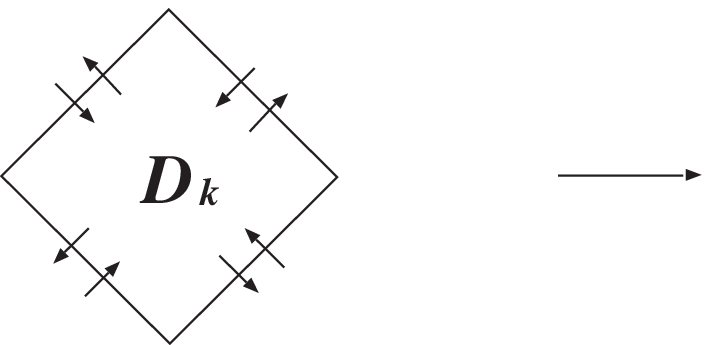}&
\includegraphics[trim=0mm 0mm 0mm 0mm, width=.3\linewidth]
{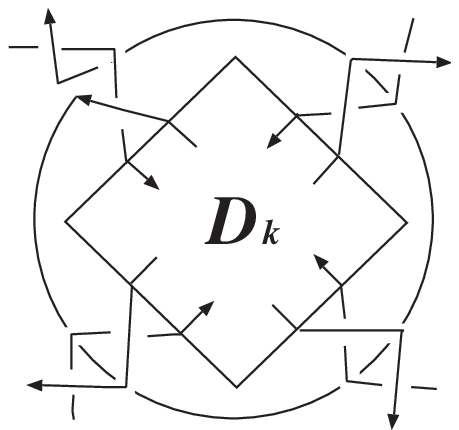}
\end{tabular}
\end{center}
\begin{center}
Fig. 4.1
\end{center}
\par
\vs
Let $F_{L_k}~~(k=1,2)$ be the Seifert surface for $L_k$~~(Fig. 4.2), and 
let ${\cal A}_{L_k}$ be the corresponding Seifert matrix of $L_k$, 
$k=1,2$.
\par
\begin{center}
\includegraphics[trim=0mm 0mm 0mm 0mm, width=.25\linewidth]
{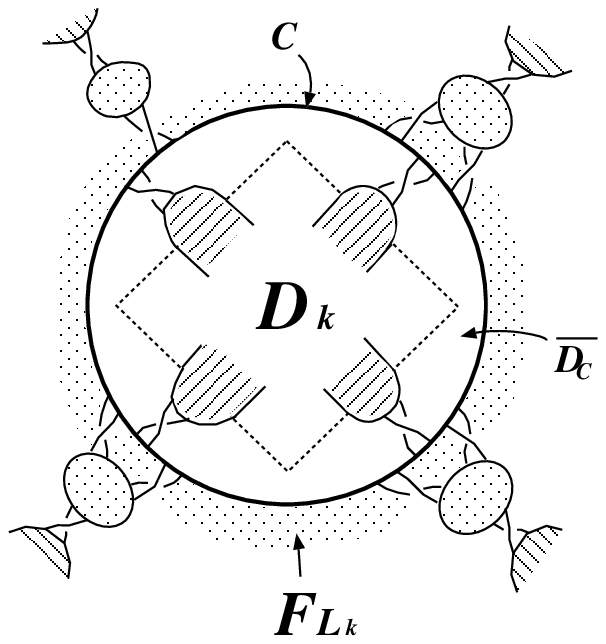}
\\
Fig. 4.2
\end{center}
\par
Let ${\xi}$ be a complex number and let 
$X_{L_k}={\xi}{\cal A}_{L_k} + {\bar{\xi}}{{\cal A}_{L_k}}^t$ 
be the Hermitian matrix that represents the Hermitian form 
$\theta(x,y)=\xi {\psi}(x,y)+ {\bar{\xi}}
{\psi}(y,x)$, $x,y \in H_1(F_{L_k};\Bbb Z)$.
\par
With the choices for Seifert surfaces $F_k$ and the bases of $H_1(F_k)$
 made above, we can formulate the main result of this section.
\par\vs
{\bf{Theorem 4.2}}~~{\em{
The characteristic polynomials of the Hermitian matrices 
$X_{L_1}$ and $X_{L_2}$ coincide.}}
\par\vs
{\bf{Proof}}~~
We consider three submodules ${\cal S}_k, {\cal R}_k$ and ${\cal
M}_k$ of $H_1(F_{L_k};\Bbb Z)$, where ${\cal S}_k, {\cal R}_k$ and ${\cal
M}_k$ are generated by the sets $X_{S_k}, X_{R_k},$ and $X_{M_k}$
of the standard generators of $H_1(F_{L_k};\Bbb Z)$ which 
live entirely in the stator part ${\bar D}$, rotor part $D_k$, 
and partially in ${\bar{D}}$ and $D_k~~(k=1,2)$, respectively.
We have the following decomposition of the module $H_1(F_{L_k};\Bbb Z)$ 
into the direct sum of its submodules, 
$H_1(F_{L_k};\Bbb Z$$)={\cal S}_k \oplus
({\cal M}_k+{\cal R}_k)~~(k=1,2)$.
\par\vs
Let $v$ denote the generator of ${\cal M}_1$
intersecting the axis $y$ of the dihedral flype $d$~~(Fig. 4.3).
There is an action of the cyclic group $\Bbb Z$$_n=<\alpha\ |\ \alpha^n=1>$ 
on ${\cal M}_1 + {\cal R}_1$ induced by the $\frac{2 \pi}{n}$-rotation 
around the
center of $D_1$. \\
The set $X_{M_1}=\{v,\alpha(v),{\alpha}^2(v), \cdots,
{\alpha}^{n-1}(v)\}$ is a generating set of ${\cal M}_1$ (not 
necessary a basis). 
%(in the case it is not a basis Theorem 4.2 follows easily, as in 
%the original Traczyk proof). 
We also identify  
${\alpha}^j(v)$ with the  generator of
${\cal M}_2$ that coincides with ${\alpha}^j(v)$ of ${\cal M}_1$ 
in ${\overline{D_C}}$.
The submodule ${\cal R}_1$ is generated by 
the set $\{{\alpha}^j(x) \mid x \in X_{R_1}, j=0,1,\cdots,n-1\}$.
Since $D_2$ is the image of $D_1$ by the dihedral flip $d$ around the
axis $y$ which crosses $v$, ${\cal R}_2$ is generated by
$\{d({\alpha}^j(x)) \mid x \in X_{R_1}, j=0,1,\cdots,n-1 \}$~~(Fig. 4.3). 
In order to compare  ${\psi}_1$ with
${\psi}_2$, we identify the generator  ${\alpha}^j(x)$ of ${\cal R}_1$ 
with  the generator 
$d({\alpha}^j(x)) \in {\cal R}_2~~(j=0,1,2,\cdots, n-1)$.     
\par

\begin{center}
\begin{tabular}{cc} 
\includegraphics[trim=0mm 0mm 0mm 0mm, width=.35\linewidth]
{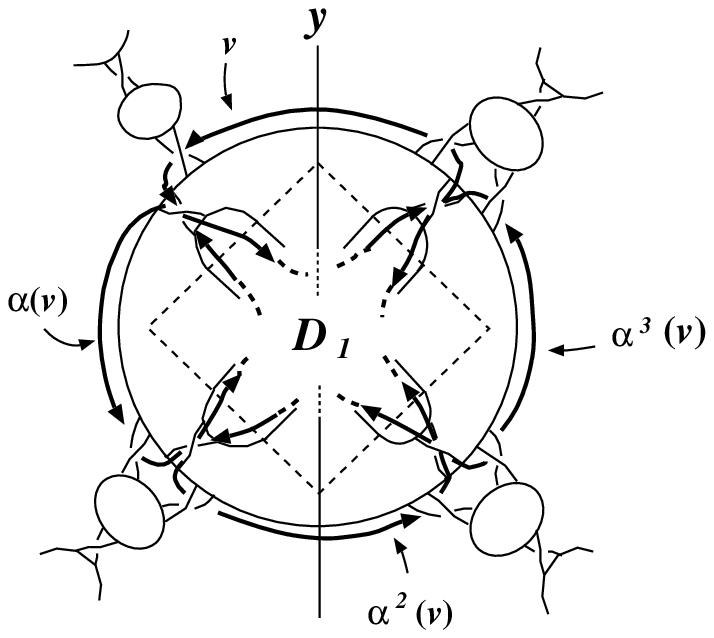}
&\hspace*{20mm}
\includegraphics[trim=0mm 0mm 0mm 0mm, width=.35\linewidth]
{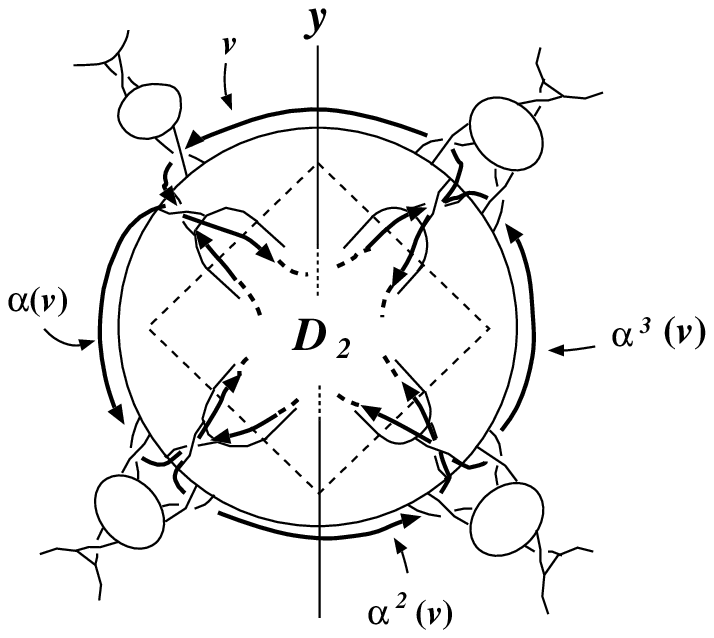}
\end{tabular}
\\
Fig. 4.3
\end{center}
\par
\vs
Using these identifications we can consider both forms 
${\psi}_1$ and ${\psi}_2$ on the same submodules 
${\cal {S,M}}$ and ${\cal R}$ (indices are no more needed) and 
derive the following relationship between them.
\par
\begin{enumerate}
\item[(1)]  ${\psi}_2(x,y)={\psi}_1(x,y)$ for all $x,y \in {\cal S}+{\cal M}$.
\item[(2)]  ${\psi}_2(x,y)={\psi}_1(y,x)$ for all $x,y \in {\cal R}$.
\item[(3)]  ${\psi}_2(x,{\alpha}^j(v))={\psi}_1({\alpha}^{-j}(v),x)$, and
\item[~~]   ${\psi}_2({\alpha}^j(v),x)={\psi}_1(x, {\alpha}^{-j}(v))$ for all $x \in {\cal R}~~(j=0,1,\cdots,n-1)$.
\item[(4)]  ${\psi}_1(x,y)={\psi}_1(y,x)=0={\psi}_2(x,y)={\psi}_2(y,x)$ for all $x
\in {\cal S},y \in {\cal R}$.
\end{enumerate}
\par\vs
Using relations (1),(2),(3),(4), we obtain the corresponding 
relations between ${\theta}_1$ and ${\theta}_2$. 
%Let ${\bf{S}}_k, {\bf{M}}_k$ and ${\bf{R}}_k~~(k=1,2)$ be the
Let ${\bf{S}}$, $ {\bf{M}}$ and ${\bf{R}}$ be the complexifications of 
subspaces $\cal S, \cal M$ and $\cal R$  of 
${\cal S} \oplus ({\cal M} + {\cal R})\otimes \Bbb C$
 respectively. 
There is a well defined involution ${\bar{}}: 
{\bf{S}} \oplus ({\bf{M}}+{\bf{R}})
\rightarrow {\bf{S}} \oplus ({\bf{M}}+{\bf{R}})$ corresponding to the
conjugation in the factor $\Bbb C$ of the tensor product. 
We denote by ${\bar{x}}$ the image of
$x \in {\bf{S}} \oplus ({\bf{M}} + {\bf{R}})$ under this involution.
The following identities follow from the identities (1)-(4) given before.
\par
\begin{enumerate}
\item[(1)]  ${\theta}_2 (x,y) =
{\theta}_1(x,y)$ for all $ x,y \in {\bf{S}} \oplus {\bf{M}}$.

\item[(2)]  ${\theta}_2 (x,y) =
{\theta}_1(y,x)={\overline{{\theta}_1(x,y)}}$ for all
generators $ x,y \in {\bf{R}}$, and 
\item[~~]  ${\theta}_2 (x,y)
={\overline{{\theta}_1({\bar{x}},{\bar{y}})}}$ for all $ x,y \in {\bf{R}}$.
\item[(3)]  ${\theta}_1(x,y)={\theta}_1({\alpha}^j(x),{\alpha}^j(y))$ for all 
generator $x,y \in {\bf{M}}+{\bf{R}}$,
\item[~~]  ${\theta}_2(x,{\alpha}^j(v))=
{\theta}_1({\alpha}^{-j}(v),x)$ for every 
generator $x$ of ${\bf{R}}$,
\item[~~]  ${\theta}_2({\alpha}^j(x),v)=
{\overline{{\theta}_1({\alpha}^{j}(x),v)}}$ for every 
generator $x$ of ${\bf{R}}$, and
\item[~~]  ${\theta}_2(x,v)=
{\theta}_2({\alpha}^j(x),{\alpha}^{-j}(v))$ for every 
generator $x \in {\bf{R}}$.
\item[(4)]  ${\theta}_k(x,y)=0$ 
for all $x \in {\bf{S}}, y \in {\bf{R}}, k=1,2$.
\end{enumerate}
\par
\vs
In order to define Hermitian matrices $H_{L_k}$ 
representing ${\theta}_k,~~(k=1,2)$,
we first choose a basis of $H_1(F_{L_k};\Bbb C)$
 that is formed using the generators of $H_1(F_{L_k};\Bbb Z)$ in the 
following way. 
Set ${\omega}_j=e^{{2 {\pi} i}{\frac{j}{n}}}~~(j=1,\cdots,n)$.
We replace the generating set $\{{\alpha}^j(v)|j=0,1,\cdots, n-1\}$~ of
${\bf{M}}$~by ~$\{{\bf{v}}_j |
{\bf{v}}_j={\displaystyle{\sum_{l=0}^{n-1}}}{\omega}_j^l{\alpha}^l(v),~j=
0,1,\cdots, n-1\}$. For ${\bf{R}}$ we consider two choices of
generating sets related by involution~ ${\bar{}}$~. We either replace 
$\{{\alpha}^j(y_p) | y_p \in X_{\cal{R}}, j=0,1,\cdots, n-1\}$~ or by
~$\{{\bf{y}}_{j,p} |
{\bf{y}}_{j,p}={\displaystyle{\sum_{l=0}^{n-1}}}{\omega}_j^l{\alpha}^l(y_p),
y_p \in X_{\cal{R}}, j=0,1,\cdots,n-1\}$~ or
~${\{\overline{{\bf{y}}_{j,p}}} |
{\overline{{\bf{y}}_{j,p}}}=
{\displaystyle{\sum_{l=0}^{n-1}}}{\overline{{\omega}_j}}^l{\alpha}^l(y_p), 
y_p \in X_{\cal{R}}, j=0,1,\cdots, n-1\}$.
\par
We obtain in this way the new generating set for $\bf{M}_k +\bf{R}_k$.
The following relationships hold: \\
\par\medskip
\begin{enumerate}
\item[(1)]  ${\theta}_k({\bf{v}}_j,{\bf{v}}_m)=0$
for $j \ne m$, where $ {\bf{v}}_j, {\bf{v}}_m \in {\bf{M}},~ k=1,2$,
\item[~~]  ${\theta}_1({\bf{x}}_{j,p},{\bf{v}}_m)={\theta}_2({\overline{{\bf{x}}_{j,p}}},{\bf{v}}_m)=0$
for $j \ne m$, where ${\bf{x}}_{j,p} \in {\bf{R}}_1, {\overline{{\bf{x}}_{j,p}}} \in {\bf{R}}_2, 
{\bf{v}}_m \in {\bf{M}}$,
\item[~~]  ${\theta}_1({\bf{x}}_{j,p},{\bf{y}}_{m,q})={\theta}_2({\overline{{\bf{x}}_{j,p}}},{\overline{{\bf{y}}_{m,q}}})$
for $j \ne m$ where ${\bf{x}}_{j,p},{\bf{y}}_{m,q} \in {\bf{R}}_1, {\overline{{\bf{x}}_{j,p}}},{\overline{{\bf{y}}_{m,q}}} \in {\bf{R}}_2$.
\item[(2)]  ${\theta}_1({\bf{x}},{\bf{y}}_{j,p})={\theta}_2({\bf{x}},{\overline{{\bf{y}}_{j,p}}})=0~~$ 
for any ${\bf{x}} 
\in {\bf{S}}, {\bf{y}}_{j,p} \in {\bf{R}}_1,{\overline{{\bf{y}}_{j,p}}} \in {\bf{R}}_2$.
\item[(3)]  ${\theta}_1({\bf{y}}_{j,p},{\bf{y}}_{j,q})=
{\overline{{\theta}_2({\overline{{\bf{y}}_{j,p}}},{\overline{{\bf{y}}_{j,q}}})}}$,~~
for any ${\bf{y}}_{j,p},{\bf{y}}_{j,q} \in {\bf{R}}_1, {\overline{{\bf{y}}_{j,p}}},{\overline{{\bf{y}}_{j,q}}} \in {\bf{R}}_2$.
\item[(4)]  ${\theta}_1({\bf{v}}_j, {\bf{y}}_{j,p})={\overline{{\theta}_2({\bf{v}}_j, {\overline{{\bf{y}}_{j,p}}})}}$ for any ${\bf{v}}_j
 \in {\bf{M}}, {\bf{y}}_{j,p} \in {\bf{R}}_1, {\overline{{\bf{y}}_{j,p}}} \in {\bf{R}}_2,$
\item[~~]  ${\theta}_1({\bf{v}}_j, {\bf{v}}_j)={\theta}_2({\bf{v}}_j,{\bf{v}}_j)$~~~for any ${\bf{v}}_j \in {\bf{M}}$.
\item[(5)]  ${\theta}_1({\bf{x}},{\bf{v}}_j)={\theta}_2({\bf{x}},{\bf{v}}_j)$ for any ${\bf{x}} \in {\bf{S}}, {\bf{v}}_j \in {\bf{M}}$.
\end{enumerate}
\par
\medskip
Take the subspace $W_j$ of ${\bf{M}} \oplus {\bf{R}}$ corresponding to
${\omega}_j$, and choose its ordered basis by taking
${\bf{v}}_j$ from ${\bf{M}}$ first\footnote{If ${\bf{v}}_j=0$, what can 
happen if the generating set $\{v,\alpha(v),{\alpha}^2(v), \cdots,
{\alpha}^{n-1}(v)\}$ is not a basis of ${\bf{M}}$, we skip 
this element when building basis of $H_1(F_L;\Bbb C)$.}  
and the rest of a basis of $W_j$ from the generating set ${\bf{y}}_{j,p}$
 of ${\bf{R}}$ 
in any order. 
%${\bf{y}}_{j,p}$ from ${\bf{R}}$ in any order. 
To obtain the ordered basis of ${\bf{M}} \oplus
{\bf{R}}$ we place the basis of $W_j$ before the basis of $W_{j+1}$ for
$j=0,1,\cdots, n-1$. Finally, we add ordered basis 
of ${\bf{S}}$. Then we have an
ordered basis of $H_1(F_L;\Bbb C)$.
We also obtain an ordered basis of $H_1(F_{L_2};\Bbb C)$ by replacing each
${\bf{y}}_{j,p}$ with ${\overline{{\bf{y}}_{j,p}}}$. 
\par
We obtain the matrices of forms ${\theta}_1$ and ${\theta}_2$ in 
ordered bases of ${\bf{S}} \oplus ({\bf{M}} + {\bf{R}})$ as described below.
\par
\begin{center}
$ H'_{L_1} = \left(
        \begin{array}{ccccc}
        B_{10} &  & {\bf{0}}  & {}^t{\overline{{S}_0}} \\
         & \ddots &  & \vdots \\
        {\bf{0}}  &  &B_{1{n-1}}  & {}^t{\overline{{S}_{n-1}}} \\
             S_0& \cdots  & S_{n-1} & S
\end{array}
       \right)$,
$ H'_{L_2} = \left(
          \begin{array}{ccccc}
        B_{20} &  & {\bf{0}}  & {}^t{\overline{{S}_0}} \\
         & \ddots &  & \vdots \\
        {\bf{0}}  &  &B_{2{n-1}}  & {}^t{\overline{{S}_{n-1}}} \\
             S_0& \cdots  & S_{n-1} & S
   \end{array}
       \right)$,
\\ 
\vspace{0.3cm}
\end{center}
In those bases, $B_{1j}$ (respectively $B_{2,j}$), where $j=0,1,\cdots,n-1$,
is the matrix of the restriction of the form ${\theta}_1$ (and 
${\theta}_2$ respectively) to the subspace $W_j$ generated by 
$\{{\bf{v}}_j\} \cup \{{\bf{y}}_{j,p} \mid y_p \in X_{{\cal
R}_1}\}$ ($\{{\bf{v}}_j \} \cup \{ {\overline{{\bf{y}}_{j,p}}} \mid y_p \in
X_{{\cal R}_1}\}$ respectively). Finally the restriction 
to the stator part, $S$, is the same for both ${\theta}_1$ and ${\theta}_2$. 
Notice that ${B_{1k}}^t=B_{2k}$, $S_l=
({\bf{s}}_{l1}$~${\bf{0}}$~$\cdots$~${\bf{0}})$, and 
 ${\bf{s}}_{l1}$ is the first column of each matrix $S_l$. 
\par
Matrices $M_k=(H'_{L_k} - {\lambda}E)~~(k=1,2)$ satisfy the
conditions\footnote{We can use Proposition 2.9, even if some vectors $w_j \in 
W_{i,j}$ may be equal to $0$. In such a case the block $W_{i,j}$ is 
orthogonal to other factors ($S$ and $W_{i,j'}$, $j'\neq j$).} 
of Traczyk's Proposition 2.9 for any real number ${\lambda}$,
\cite{T}.
Thus $\det(M_1)=\det(M_2)$ for any real $\lambda$. So the determinants
are equal for any complex number $\lambda$ as well.
\par
\hfill$\Box$

\section{Counterexamples}
%\begin{section{Counterexample}}
It was proven in \cite{APR} that any pair of oriented $3$- 
or $4$-rotant links share the
 same Homflypt polynomials (in particular, Conway polynomials).
In \cite{T} Traczyk showed that a pair of orientation-preserving $n$-rotant 
links share the same Conway polynomial.
On the other hand, 
for orientation-reversing $n$-rotants $(n \geq 6)$, the invariance 
was still an open question. 
We present, in this section, an example of a pair of $6$-rotant knots 
with different Conway polynomials and different Jones
 polynomials. Therefore, the invariance in \cite{APR} of Conway
 polynomial and  the Jones polynomial for the orientation-reversing rotant
 links is the best possible. We should also stress that rotants described 
in Fig. 5.1 
have different Jones and Conway polynomials, however they share 
the same determinant and the same
 homology of the corresponding double branched covers. 
\par\vs
\begin{center}
\begin{tabular}{cc} 
\includegraphics[trim=0mm 0mm 0mm 0mm, width=.3\linewidth]
{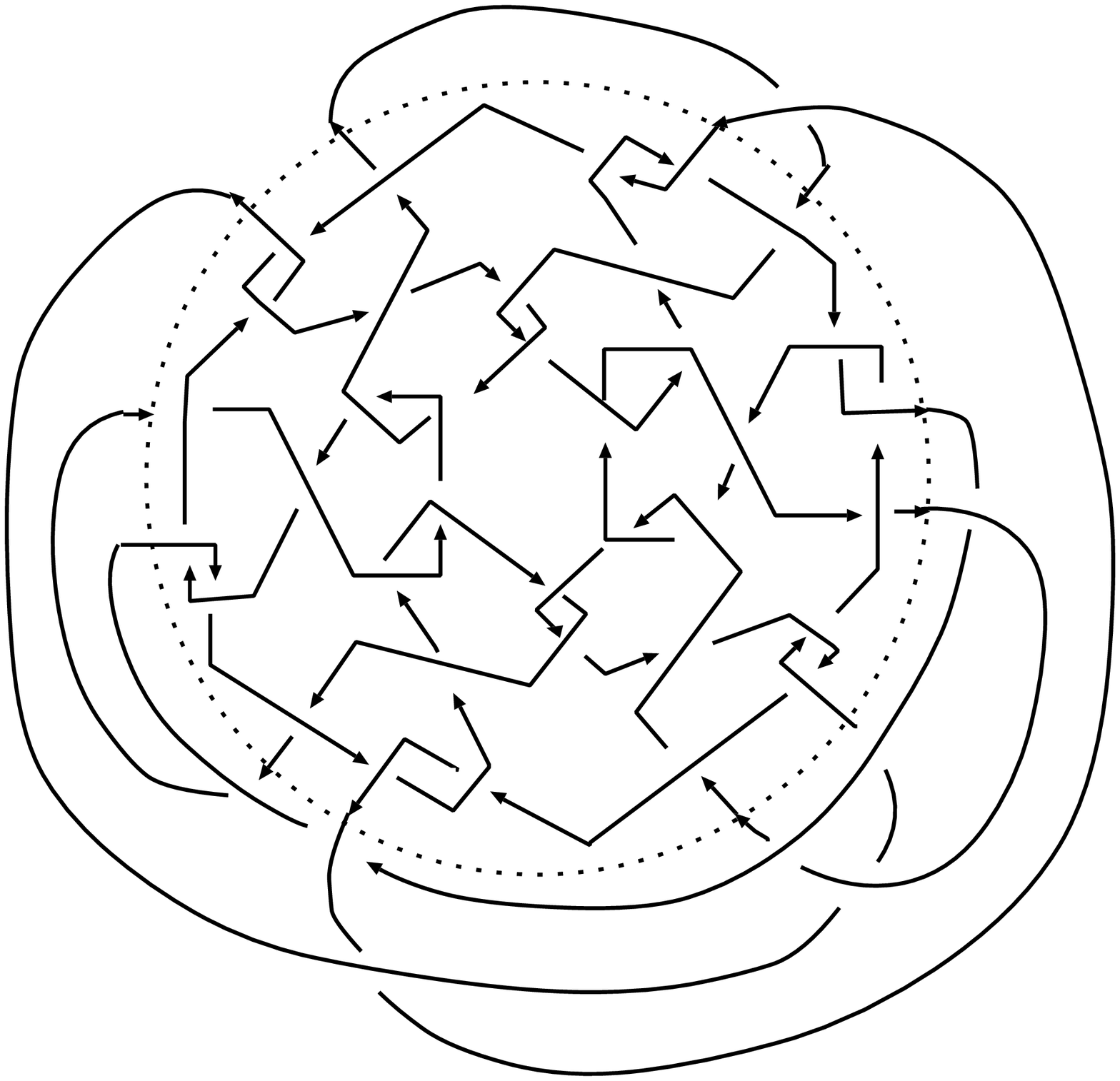}
&\hspace*{10mm}
\includegraphics[trim=0mm 0mm 0mm 0mm, width=.3\linewidth]
{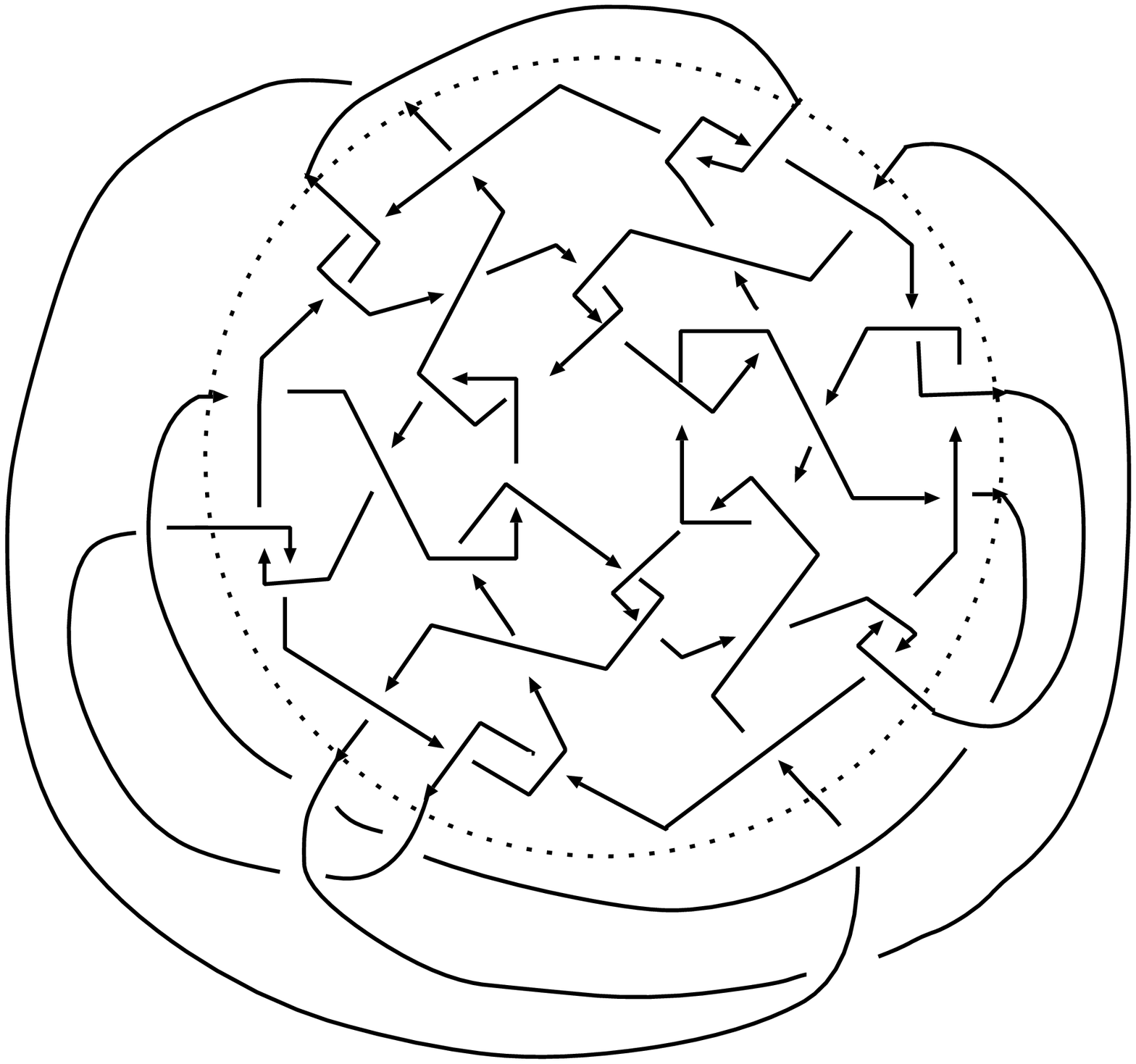}\\
$L_1$&\hspace*{10mm}$L_2$
\end{tabular}
\par
Fig.5.1
\end{center}
\par\vs
Let $L_1$ and $L_2$ be the knots (6-rotants) illustrated in Fig. 5.1.
Using program KNOT \cite{K}, we have the following.
\par\vs
Conway polynomials (with the skein relation 
${\bigtriangledown}_{L_+}-{\bigtriangledown}_{L_-}=z{\bigtriangledown}_{L_0}$) 
are different:
\par
${\bigtriangledown}_{L_1}(z)=1+3z^2-37z^4+17z^6-3z^8-2z^{10}-59z^{12}-34z^{14}-55z^{16}-48z^{18}-10z^{20}$
\par
~~~~~~~~~~~~~~$-4z^{22}-z^{24}$, 
\par\vs\vs
and
\par\vs\vs
${\bigtriangledown}_{L_2}(z)=1+3z^2-25z^4-116z^8-57z^{10}-174z^{12}-157z^{14}-119z^{16}-102z^{18}$
\par
~~~~~~~~~~~~~~$-37z^{20}-8z^{22}-z^{24}$.

\par\vs\vs\par\vs\vs
Jones polynomials (with the skein relation  
$t^{-1}V_{L_+}-tV_{L_-}=(\sqrt{t}-1/{\sqrt{t}})V_{L_0}$) are different:
\par\vs
$V_{L_1}=t^{23}-16t^{22}+131t^{21}-713t^{20}+2881t^{19}-9193t^{18}+24058t^{17}-52926t^{16}$
\par
~~~~~~~~$
+99534t^{15}-161854t^{14}+229195t^{13}-283357t^{12}+304679t^{11}-280476t^{10}
$
\par
~~~~~~~~$+211413t^{9}-112418t^{8}+7697t^{7}+77824t^{6}-127092t^{5}+136195t^{4}-114114t^{3}$
\par
~~~~~~~~~$+77214t^{2}-41391t+16087-2934t^{-1}-1501t^{-2}+1760t^{-3}-954t^{-4}$
\par
~~~~~~~~~$+343t^{-5}-84t^{-6}+13t^{-7}-t^{-8}$, 
\par\vs\vs
and
\par\vs\vs
$V_{L_2}=t^{23}-16t^{22}+131t^{21}-713t^{20}+2881t^{19}-9193t^{18}+
24057t^{17}-52919t^{16}$
\par
~~~~~~~~~$+99503t^{15}-161752t^{14}+228932t^{13}-282808t^{12}
+303730t^{11}-279098t^{10}$
\par
~~~~~~~~~$+209727t^{9}-110701t^{8}+6314t^{7}+78540t^{6}-126958t^{5}+
135242t^{4}$
\par
~~~~~~~~~$-112578t^{3}+75451t^{2}-39756t+14823-2118t^{-1}-1933t^{-2}+
1941t^{-3}$
\par
~~~~~~~~~$-1010t^{-4}+354t^{-5}-85t^{-6}+13t^{-7}-t^{-8}$.
\par\vs\vs
Their homology groups are the same:
$H_1(M^2_{L_1}; {\Bbb Z})$$=H_1(M^2_{L_2}; {\Bbb Z})$
$={\Bbb Z}$$/3 \oplus {\Bbb Z}$$/397449$. Their determinants
%and signatures
 coincide as well: 
${\Delta}_{L_1}(-1)={\Delta}_{L_2}(-1)=-1192347$ (here ${\Delta}_{L}(t)=
\nabla_L(z)$ for $z=\sqrt{t}-\frac{1}{\sqrt{t}}$.).
%$\sigma(L_1)= \sigma(L_2) = 4k+2$.

%\input{proof-th2.tex}
%\input{proofclaim.tex} 
%\par\vs\vs
%{\bf{Claim}}: ~$e(F_{L_1})=e(F_{L_2})$
%\par\vs
%{\it{Proof of Claim}}
%\par
\par\vs
%%%%%%%%%% References
\end{document}